\documentclass[a4paper,notitlepage,twoside,reqno,11pt]{amsart}

\usepackage{mathrsfs}
\usepackage{yhmath}
\usepackage{tikz}
\usepackage[hmargin=1in,vmargin=1in]{geometry}
\usepackage{amsmath}
\usepackage{amssymb}
\usepackage{amsfonts}
\usepackage{amsthm}
\usepackage{bbm,pifont,latexsym}
\usepackage{anysize}  
\marginsize{3cm}{3cm}{3cm}{3cm}

\usepackage{dcolumn,indentfirst,color,hyperref}
\usepackage{amsmath,amssymb,amscd,amsthm,amsfonts,mathrsfs}
\usepackage{color,graphicx,xcolor,graphics}
\usepackage{titlesec} 
\usepackage{titletoc} 
\titlecontents{section}[20pt]{\addvspace{2pt}\filright}
              {\contentspush{\thecontentslabel. }}
              {}{\titlerule*[8pt]{}\contentspage}

\usepackage{fancyhdr} 
\newtheorem{thm}{Theorem}[section]

\newtheorem{lem}[thm]{Lemma}
\newtheorem{cor}[thm]{Corollary}
\newtheorem{prop}[thm]{Proposition}

\newtheorem{rmk}[thm]{Remark}
\newtheorem{defi}{Definition}[section]

\newcommand{\wt}{\widetilde}
\newcommand{\mb}{\mathbb}
\newcommand{\mc}{\mathcal}

\newcommand{\tu}{\textup}
\newcommand{\ol}{\overline}
\newcommand{\es}{\emptyset}

\makeatletter\@addtoreset{equation}{section}\makeatother 
\titleformat{\section}{\centering\normalsize}{\textsc{\thesection.}}{1em}{\textsc}
\titleformat{\subsection}{\normalsize}{\thesubsection.}{1em}{\textbf}

\title
 {Criterion for rays landing together}
\begin{document}

\author{Jinsong Zeng}
\address{Jinsong Zeng, School of Mathematics and Information Science, Guangzhou University, Guangzhou 510006, P. R. China}
\email{jinsongzeng@163.com}

\begin{abstract}
We give a criterion to determine when two external rays land at the same point for polynomials with locally connected Julia sets. As an application, we provide an elementary proof of the monotonicity of the core entropy  along arbitrary veins of the Mandelbrot set.
\end{abstract}



\maketitle

\section{Introduction}
Let $f:\mb{C}\to\mb{C}$ be a polynomial of degree $d\geq 2$. The \emph{basin of infinity} $\Omega_f$ consists of the points whose orbits escape to infinity under the iterations of $f$. The \emph{Julia set} $J_f$ is the boundary of $\Omega_f$. One of the main interests in complex dynamics is to study the combinatorics and topology of the Julia set as well as their relationship with the dynamics of $f$ on the Julia set.

If the Julia set $J_f$ is connected, it is known from \cite[Theorem 9.5]{Mi06} that there is a unique conformal isomorphism $\Psi_f:\Omega_f\to\mathbb{C}\setminus\ol{\mb{D}}$ satisfying the properties: $\Psi_f(f(z))=\Psi_f(z)^d$ for all $z\in \Omega_f$ and $\tu{lim}_{z\to\infty}\Psi_f(z)/z=1$. The map $\Psi_f$ is called the \emph{B\"{o}ttcher map} of $f$. The \emph{(external) ray} $R(\theta)$ of an angle $\theta\in\mb{T}:=\mb{R}/\mb{Z}$ is defined by
$$R(\theta):=\Psi^{-1}_f(\{re^{2\pi i\theta}:r>1\}).$$

If the Julia set $J_f$ is connected and locally connected, by Carath\'{e}odory's theory, the continuous extension of the inverse $\Psi_f^{-1}$ over the boundary $\partial \mb{D}$ onto the whole Julia set exists and is a semi-conjugacy between the maps $z^d:\partial \mb{D}\to \partial\mb{D}$ and $f:J_f\to J_f$. In this case, each ray $R(\theta)$ \emph{lands} at a point in the Julia set, i.e., the limit $\tu{lim}_{r\to 1^+}\Psi^{-1}_f(re^{2\pi i\theta})$ exists.

Two questions arise naturally:


\begin{enumerate}
	\item[(Q1)] Given a point in the Julia set, how many rays are there landing at this point?
	\item[(Q2)] Given two rays $R(\theta_1)$ and $R(\theta_2)$, under what conditions do they land at the same point?
\end{enumerate}

A point is called \emph{wandering} if its forward orbit under the iteration of $f$ is infinite.
In \cite[Theorem 1.1]{Ki02}, it was shown that the number of rays landing at a wandering point in the Julia set is bounded by $2^d$. By working on Thurston's \emph{invariant lamination} \cite{Th85} and  Levin's \emph{growing tree} \cite{Le98} (a generalization of Hubbard tree), Blokh and Levin proved the following general result \cite[Theorem A]{BL02}.

\begin{thm}\label{thm:1}
Let $f$ be a polynomial of degree $d\geq 2$ with locally connected Julia set.
Let $z_1,\cdots,z_m$ be wandering points in the Julia set satisfying that their forward orbits are pairwise disjoint and avoid the critical points of $f$. Then
\begin{equation}\label{estimate}
\sum_{1\leq i\leq m}(v(z_i)-2)\leq d-2,
\end{equation}
where $v(z_i)$ is the number of rays landing at $z_i$.
\end{thm}
A point in $\mb{C}$ is \emph{branched}, if it is in the Julia set and is the landing point of more than two rays. Thurston proved that there is no wandering branched point for any quadratic polynomial \cite{Th85} (it also follows directly from \eqref{estimate}). However, this is not true in general. Blokh and Oversteegen constructed cubic polynomials whose Julia sets contain wandering branched points; see \cite{BO08}.

To study question (Q2), we need the notion of \emph{critical portrait}, a collection of subsets of the unit circle $\mb{T}$; see Definition \ref{def:critical_portrait} or \cite{BFH92,Po09} for more details.
A critical portrait naturally divides $\mb{T}$ into $d$ pieces. Associated to these pieces, all but countably many angles in $\mathbb{T}$ have well-defined itineraries under the iterations of $m_d:\theta\mapsto d\,\theta\tu{ mod }\mb{Z}$ (Definition \ref{itinerary}).

For a polynomial with all critical points strictly preperiodic, if two angles have the same itinerary, then the rays corresponding to these two angles land at the same point; see \cite[Section 3]{BFH92}. This result was extended by Poirier to postcritically finite polynomials; see \cite[Corollary 5.9]{Po09}. Both of their proofs make essential use of the fact that $f:J_f\to J_f$ is expanding with respect to the \emph{orbifold metric}.

For a polynomial with all cycles repelling and with connected Julia set, if it is \emph{visible} meaning that it can be approached by polynomials in the shift locus, then the above criterion for two rays landing at the same point holds again (deduced from \cite[Theorem 1]{Ki05}). This result is due to the combinatorial continuity of complex polynomial dynamics. 

\subsection{Our results}
We are mainly concerned with the questions (Q1) and (Q2).
Our first result is to give a brief and new proof of Theorem \ref{thm:1}. The method employed here is quite different from that of Blokh and Levin. It is based on the combinatorial analysis of the  \emph{orbit portraits} of wandering branched points. The notion of orbit portraits was introduced by Goldberg and Milnor in \cite{GM93}. The second and main result of this paper is as follows.
\begin{thm}[Main Theorem]\label{thm:main_thm}
Let $f$ be a polynomial of degree $d\geq 2$ with locally connected Julia set. If two angles $\theta_1$ and $\theta_2$ have the same itinerary with respect to a critical portrait, then the landing points of $R(\theta_1)$ and $R(\theta_2)$ either coincide or belong to the boundary of a Fatou domain, which is eventually iterated onto a Siegel disk.
\end{thm}

Note that Zakeri \cite[Theorem 5]{Za00} proved that, for a \emph{Siegel quadratic polynomial}, i.e. $f:z\mapsto z^2+c$ with a fixed Siegel disk, there is no branched point; and if two rays land at $z$, then $z$ eventually hits the unique critical point.

As a consequence of Theorem \ref{thm:main_thm}, we have
\begin{cor}[No wandering continuum in $J_f$]\label{cor:no_wandering}
Let $f$ be a polynomial of degree $d\geq 2$ whose Julia set $J_f$ is locally connected. Then any continuum $C$ in $J_f$ must have $f^{m}(C)\cap f^{m+n}(C)\neq\es$ for some $m\geq0$ and $n\geq 1$.
\end{cor}
By a \emph{continuum} we mean a compact, connected and non-singleton subset of $\mb{C}$.
Blokh and Levin also proved Corollary \ref{cor:no_wandering} in \cite[Theorem C]{BL02}. By using the Yoccoz puzzle technique, Kiwi showed that, when $f$ has no irrational neutral periodic cycle, $J_f$ is locally connected if and only if $f$ has no wandering continua in $J_f$; see \cite[Theorem 5.12]{Ki04}. Note that Yoccoz puzzles always come from $f$-invariant graphs, while in our proof the graphs induced by critical portraits are certainly not $f$-invariant. 


The main motivation for answering the questions (Q1) and (Q2) is to study the core entropy of polynomials, which was first introduced and explored by Thurston \cite{Th14}.
For a postcritically finite polynomial $f$, the \emph{core entropy} $h(f)$ is defined as the entropy of the action of $f$ on its \emph{Hubbard tree}. By definition, the Hubbard tree is the convex hull of the critical orbits within the \emph{filled Julia set}, i.e., the complement of the basion of infinity. 

The \emph{biaccessible set} $\tu{Acc}(f)$ is the set of all \emph{biaccessible angles}, which are the angles $\theta\in\mb{T}$ such that both $R(\theta)$ and $R(\theta')$ land at the same point for some $\theta'\neq \theta$.
In the case that $f$ is postcritically finite, the core entropy $h(f)$ is closely related to the Hausdorff dimension of $\tu{Acc}(f)$ in the following way
\begin{equation}\notag
h(f)=\tu{log}\,d\cdot\tu{H.dim}\, \tu{Acc}(f);
\end{equation}
see \tu{\cite[Theorem 7.1]{Ti13}}. For a general polynomial with locally connected Julia set, it may occur that the Hubbard tree is infinite, but
the core entropy $h(f)$ can be defined as $\tu{log }d\cdot \tu{H.dim}\, \tu{Acc}(f)$ (cf. \cite{BS},\cite[Appendix A]{DS},\cite{Ti13}).

As an application of Theorem \ref{thm:main_thm}, in the last section we provide an alternative (elementary) proof of the following theorem about the monotonicity of the core entropy for the family of quadratic polynomials; see \cite{DS,Ju13} for other proofs.



\begin{thm}\label{thm:mon}
For all $f_c$ and $f_{c'}$ in the family 
$$\mc{F}=\{f_c:z\mapsto z^2+c\tu{ having locally connected Julia set without a Siegel disk}\},$$ if $f_c\prec f_{c'}$, then $\tu{Acc}(f_c)\subseteq \tu{Acc}(f_{c'})$ and so $h(f_c)\leq h(f_{c'})$.
\end{thm}

Here we say $f_c\prec f_{c'}$ if $I_c\supseteq I_{c'}$, where $I_c$ is the \emph{characteristic arc} of $f_c$; see Definition \ref{def:characteristic} or \cite[Lemma 2.6]{Mi00}. The family $\mc{F}$ contains all \emph{veins} within the Mandelbrot set; and in the sense of the partial relation $``\prec"$ on $\mc{F}$, all elements in a vein form a totally ordered set \cite{Do93,PR,Ri00}. Therefore, the monotonicity of the core entropy along arbitrary veins is established. 
One may refer to \cite{Do95,DS,Ju13,Ti14} for the continuity of core entropy in the Mandelbrot set.

\subsection{Outline of the paper}
In Section \ref{Wandering}, we give a brief and new proof for Theorem \ref{thm:1} by studying the dynamics of portraits under the sector maps. Regulated arcs for topological polynomials are studied in Section \ref{regulated_arcs}, where we establish the rigidity (existence and uniqueness) of the arcs and the density of branched points in certain regulated arcs, namely the clean arcs. In Section \ref{critical_portraits} we study the properties of critical portraits, which play an essential role in the proof of Theorem \ref{thm:main_thm}. With all the preparations in Section \ref{Wandering}, \ref{regulated_arcs} and \ref{critical_portraits}, we prove Theorem \ref{thm:main_thm} in Section \ref{proof_main_thm}, without using Thurston's invariant lamination. As an application of Theorem \ref{thm:main_thm}, we give an elementary proof for Theorem \ref{thm:mon} in Section \ref{application}.

\vskip 1cm

\emph{Acknowledgment. }
The main results of this paper are from a revised version of the unpublished part of the author's Ph.D thesis.
We would like to thank Weiyuan Qiu and Lei Tan (1963-2016), the thesis advisors of the author, for their many helpful discussions on this paper. We also thank the anonymous referees for useful comments and careful reading of this paper. The author was partially supported by the  China Scholarship Council and the NSFC under grant No.11801106 during the preparation of this paper.$\\$

\section{Wandering Orbit Portrait}\label{Wandering}
In this section we establish that, for wandering branched points with pairwise disjoint forward orbits, the total number of rays landing at them has an upper bound in terms of the degree $d$, as stated in Theorem \ref{thm:1}. 

We always assume in this section that $f$ is a polynomial of degree $d\geq 2$ with locally connected Julia set. A point in the Julia set with at least two landing rays is called \emph{biaccessible}, and is called \emph{branched} if there are more than two rays landing at it.

\begin{defi}[Portrait]
Let $z$ be a biaccessible point. A finite subset $T$ of the unit circle $\mb{T}:=\mb{R}/\mb{Z}$ with $\#T\geq 2$ is called a \emph{portrait} of $z$ if for each $\theta\in T$, the ray $R(\theta)$ lands at $z$.

Let $T$ be a portrait of $z$. The set $\cup_{\theta\in T}R(\theta)$ cuts $\mb{C}\setminus\{z\}$ into $\#T$ parts. These parts are called the \emph{sectors} of the portrait $T$ based at $z$. Associated to each sector $S$, there is an open interval $I$ of $\mb{T}$ formed by 
$$I=\{\theta\in\mb{T}: R(\theta)\subseteq S\}.$$
The \emph{angular size} $l(S)$ is defined as the length of the interval $I$.
\end{defi}


\begin{lem}[{Sectors based at distinct points}]\label{lem:distinct_portrait}
Let $T,T'$ be two portraits of distinct points $z$ and $z'$. Let $S$ (resp. $S'$) be the sector of $T$ (resp. $T'$) containing $z'$ (resp. $z$). Then it holds that
$$Q\subseteq S\tu{ and }l(Q)<l(S),$$
where $Q$ runs over all sectors of $T'$ except $S'$.
\end{lem}

\begin{proof}
	It is clear that $\partial S'\subseteq S$. The domain $W:=\mb{C}\setminus\ol{S'}$ is a subset of $S$. Since each sector $Q$ of $T'$ except $S'$ belongs to $W$, it follows that $Q\subsetneq W\subsetneq S$. Hence $l(Q)<l(S)$.
\end{proof}

The image of a sector under $f$ may not be a sector. However, locally around a non-critical point $z$ in the Julia set, the action of $f$ induces a bijection, called a \emph{sector map}, between sectors based at $z$ and sectors based at $f(z)$. The precise definition is as follows.

\begin{defi}[Sector map]
	Let $z$ be a non-critical point in the Julia set and $T$ be a portrait of $z$. Since $f$ is homeomorphic on a neighborhood $V$ of $z$ and carries the rays landing at $z$ to those at $f(z)$, it holds that $T':=m_d(T)$ is a portrait of $w:=f(z)$, where $m_d: \mb{T}\to \mb{T}$ sends each $\theta$ to $d\theta\tu{ mod }\mb{Z}$. Moreover, we have $$f(\cup_{\theta\in T}R(\theta))=\cup_{\theta'\in T'}R(\theta').$$
	For each sector $S$ of $T$, the image $f(S\cap V)$ belongs to a unique sector, say $\tau(S)$, of $T'$. We then call $\tau$ \emph{the (local) sector map} at $z$.
\end{defi}

In the following, we include a lemma due to Goldberg and Milnor \cite{GM93}, which describes the basic relations between the number of critical points and values within a sector, angular sizes and sector maps.

\begin{lem}[Properties of sector maps]\label{lem:portrait_map}
Let $T$ be a portrait of a non-critical point $z$ in the Julia set. Let $S$ be a sector of $T$.
Then
\begin{enumerate}
	\item The sector map $\tau$ is a bijection  from all sectors of $T$ to those of $T':=m_d(T)$.
	
	\item $l(\tau(S))= dl(S)\tu{ mod }\mb{Z}$. The integer $n_0:=dl(S)-l(\tau(S))$ is the number of critical points, counting multiplicity, of $f$ contained in $S$.
	
	\item If $n_0\geq 1$ or $l(\tau(S))\leq l(S)$, then $\tau(S)$ contains at least one critical value.
	
	\item If $l(S)<1/d$, then the restriction $f:S\to \tau(S)$ is a homeomorphism.
\end{enumerate}
\end{lem}

\subsection{Dynamics of portraits}
A point in the Julia set is called \emph{wandering} if its forward orbit is infinite. Let $T_0$ be a portrait of a wandering and non-critical branched point $z$ with $u:=\#T_0\geq 3$. We assume that the forward orbit of $z$ avoids the critical points of $f$. Then we get a sequence of portraits $T_n:=m_d^n(T)$ of $f^n(z)$. 

We denote by $S_{1,n},\cdots,S_{u,n}$ the $u$ sectors of $T_n$ for each $n\geq 0$, such that
\begin{equation}\label{eq:less}
l(S_{1,n})\geq l(S_{2,n})\geq \cdots\geq l(S_{u,n}).
\end{equation}
\begin{lem}\label{lem:sector_converging_to_zero} For $3\leq k\leq u$, we have
	$\lim_{n\to\infty}l(S_{k,n})=0.$
\end{lem}
\begin{proof}
It suffices to show that $\lim_{n\to\infty}l(S_{3,n})=0$ by \eqref{eq:less}. We assume this is not true, and then there exists a subsequence $l(S_{3,n_i})$ of angular sizes with limit $a>0$ such that
	$$\frac{5}{6}a\leq l(S_{3,n_i})\leq\frac{7}{6}a,$$
for each $i\geq 0$.
	Therefore at most finitely many of sectors among $\{S_{3,n_i}\}_i$ are pairwise disjoint. Thus we can select two nested sectors, say $S_{3,n_0}$ and $S_{3,n_1}$, such that 
	$$S_{3,n_0}\supsetneq S_{3,n_1}\tu{ and }f^{n_1}(z)\in S_{3,n_0}.$$
	According to Lemma \ref{lem:distinct_portrait}, at least one of $S_{1,n_1}$ and $S_{2,n_1}$ is contained in $S_{3,n_0}$. Therefore
	$$l(S_{3,n_0})\geq l(S_{3,n_1})+\tu{min}\,\{l(S_{1,n_1}),l(S_{2,n_1})\}\geq \frac{5}{6}a+\frac{5}{6}a.$$
	This contradicts the assumption that $l(S_{3,n_0})\leq \frac{7}{6}a$. 
\end{proof}
A sector containing critical points (resp. critical values) is called a \emph{critical sector} (resp. \emph{critical-value sector}).
The image of a critical sector under the sector map must be a critical-value sector according to Lemma \ref{lem:portrait_map}. 

By the expansion of $m_d$ on $\mb{T}$,
any sector will eventually be iterated to a critical-value sector by the sector maps. 
The main focus lies in the very moment when a ``wide" sector becomes ``narrow". From Lemma \ref{lem:portrait_map}, this phenomenon happens only when a critical sector is sent to a critical-value sector.

We adopt the notations $z, T, T_n, u$ and $S_{k,n}$ as above. For a small $\epsilon>0$ and an integer $3\leq k\leq u$, set
$$N_{\epsilon,k}(T):=\tu{min}\,\{n:l(S_{k,n})<\epsilon\}.$$
Such an integer $N_{\epsilon,k}(T)$ always exists according to Lemma \ref{lem:sector_converging_to_zero}.  By analyzing the changes of the angular sizes from $T_{N-1}$ to $T_{N}$ for $N=N_{\epsilon,k}(T)$, we have the following crucial lemma.

\begin{lem}\label{lem:critical_value_sector}
	Let $\epsilon, k$ and $N$ be given as above. We consider the sector map  $$\tau:\{S_{1,N-1},\cdots,S_{u,N-1}\}\to \{S_{1,N},\cdots,S_{u,N}\}.$$
	It holds that
	\begin{enumerate}
		\item $l(S_{k-1,N})>\epsilon$;
		\item there exists a critical-value sector among $S_{k,N},\cdots, S_{u,N}$;
		\item $N'\neq N$ for another $3\leq k'\leq u$ and $N':=N_{\epsilon,k'}(T)$.
	\end{enumerate} 
\end{lem}
\begin{proof}By Lemma \ref{lem:sector_converging_to_zero},
	choose $\epsilon_0$ so small that, for any $\epsilon<\epsilon_0$ and $n\geq N-1$, the angular sizes of $S_{3,n},\cdots, S_{u,n}$ are less than $\frac{1}{4ud}.$ Thus for $3\leq i\leq u$, $f$ maps $S_{i,N-1}$ conformally onto the sectors of $T_N$. 
	We see that $$\tau\{S_{3,N-1},\cdots,S_{u,N-1}\}\subseteq\{S_{2,N},\cdots,S_{u,N}\}.$$
	It is because
	$$l(S_{1,N})\geq \frac{1}{2}(l(S_{1,N})+l(S_{2,N}))=\frac{1}{2}(1-\sum_{3\leq i\leq u}l(S_{i,N}))\geq \frac{1}{2}(1-\frac{1}{4d})>\frac{3}{8}>l(\tau(S_{j,N-1}))$$
	for all $3\leq j\leq u$. Then one of $S_{1,N-1}$ and $S_{2,N-1}$, say $S_{2,N-1}$, is sent into $\{S_{2,N},\cdots,S_{u,N}\}$. In what follows, one may assume without loss of generality that none of the equalities in \eqref{eq:less} holds for $S_{1, N-1}, \cdots, S_{u,N-1}$.
	
	We claim that $S_{2,N-1}$ is critical. Otherwise, the sector map $\tau$ would preserve the order of the angular sizes of $S_{i,N-1}$ for $2\leq i\leq u$. Then $l(S_{k,N})=d\cdot l(S_{k,N-1})\geq \epsilon$, which is a contradiction. 

	
	Next, we have $\tau(S_{2,N-1})\in\{S_{k,N},\cdots,S_{u,N}\}$. In fact, if $\tau(S_{2,N-1})\in\{S_{2,N},\cdots,S_{k-1,N}\}$, then the appearance of the critical-value sector $\tau(S_{2,N-1})$ does not break the order of the smallest $u-k+1$ sectors $S_{k,N-1},\cdots,S_{u,N-1}$ under the action of $\tau$. Therefore one has $l(S_{k,N})=d\cdot l(S_{k,N-1})\geq \epsilon$, which is impossible. Hence (2) follows.
	
	The discussions above also indicate that $l(S_{k-1,N})=l(\tau(S_{k,N-1}))\geq d\epsilon>\epsilon$. This gives (1).
	Moreover, it holds that $S_{i,N}=\tau(S_{i,N-1})$ for all $k+1\leq i\leq u$. Hence (3) follows.
\end{proof}
\begin{prop}\label{prop:key_inequality}
	Let $T^{(1)},\cdots,T^{(m)}$ be portraits of wandering branched points $z_1,\cdots, z_m$ with $u_i:=\#T^{(i)}\geq 3$. Suppose that the forward orbits of $z_1,\cdots,z_m$ are pairwise disjoint and avoid the critical points of $f$. Then
	$$\sum_{1\leq i\leq m}(u_i-2)\leq d-2.$$
\end{prop}
\begin{proof}
	Let $\epsilon$ be small enough satisfying Lemma \ref{lem:critical_value_sector} for each $T^{(i)}$. Let $N_{i,k}:=N_{\epsilon,k}(T^{(i)})$ and $V_{i,k}$ be the critical-value sectors based at $f^{N_{i,k}}(z_i)$ in Lemma \ref{lem:critical_value_sector} (2) for all
	$$1\leq i\leq m\tu{ and } 3\leq k\leq u_i.$$ According to Lemma \ref{lem:critical_value_sector} (3), the $s:=\sum_{1\leq i\leq m}(u_i-2)$ points $\{f^{N_{i,k}}(z_i)\}$ are all distinct.
	
	Since the total number of the critical values of $f$ is at most $d-2$, it suffices to show that the $s$ sectors in $\{V_{i,k}\}$ are pairwise disjoint. As sectors are either nested or disjoint, we are left to prove that there is no sector contained in another one among $\{V_{i,k}\}$.
	
	 For otherwise, suppose $V_{i,k}\subseteq V_{i',k'}$ with $(i,k)\neq (i',k')$. Let us consider the image of $T^{(i)}$ under the map $m_d^{N_{i,k}}$. We denote it by $T$.  Let $A_1$ and $A_2$ be the two sectors of $T$ with the largest angular sizes. According to Lemma \ref{lem:critical_value_sector} (2), we have $V_{i,k}\neq A_1,A_2$. 
	 Lemma \ref{lem:distinct_portrait} implies that one of $A_1$ and $A_2$ is contained in $V_{i',k'}$. Then we have
	 $$l(V_{i',k'})\geq\tu{min}\,\{l(A_1),l(A_2)\}>\epsilon,$$
	 where the latter inequality is from Lemma \ref{lem:critical_value_sector} (1). This is a contradiction and we complete the proof.
\end{proof}
\noindent{\textit{Proof of Theorem \ref{thm:1}}}. The theorem follows immediately from Proposition \ref{prop:key_inequality}.
\hfill $\square$

\begin{cor}\label{cor:finite_branched}
	Let $f$ be a polynomial of degree $d \geq 2$ with locally connected Julia set. Then there are at most $d-2$ wandering branched points in the Julia set with disjoint forward orbits.
\end{cor}
This corollary will be needed in Section \ref{proof_main_thm}.

\section{Regulated arcs}\label{regulated_arcs}

This section is devoted to the generalization of the notion of \emph{regulated arcs}, introduced by Douady-Hubbard \cite{DH84} for postcritically finite polynomials, to all polynomials with locally connected Julia sets.
  



\begin{lem}\label{lem:local_connectivity}
	Let $f$ be a polynomial of degree $d\geq 2$ with locally connected Julia set. Then
	\begin{enumerate}
		\item any bounded Fatou domain is a Jordan domain;
		\item the number of Fatou domains with diameters greater than a given positive number $\epsilon_0$ is at most finite.
	\end{enumerate}
\end{lem}
\begin{proof}
	By the property of local connectivity, $(2)$ holds. Moreover, the boundary of a bounded Fatou domain $U$ is locally connected. If $\partial U$ is not a Jordan curve, then there is a Jordan curve $\gamma$ and a point $p\in\partial U$ such that $\gamma\cap \partial U=\{p\}$ and $\gamma\setminus\{p\}\subseteq U$ with both components of $\mb{C}\setminus\gamma$ intersecting the Julia set. This is impossible, as the Julia set is the boundary of the basin of infinity.
\end{proof}

\begin{defi}
Let $f$ be a polynomial with locally connected Julia set. For any bounded Fatou domain $U$, choose a point $c_U$ in $U$ as the \emph{center} of $U$, and a homeomorphism $\phi_U:\ol{U}\to \ol{\mb{D}}$ that sends $c_U$ to $0$. The pullbacks of radial rays, i.e., $R_U(\theta):=\phi_U^{-1}(e^{2\pi i\theta}[0,1))$, are called the \emph{internal rays} of $U$.
A \emph{closed internal ray} means the union of this ray and its landing point.
The collection $\Phi_f:=\{(c_U,\phi_U)\}_U$ is a \emph{parameterization} of bounded Fatou domains (of $f$).	
\end{defi}

Since the image of an internal ray under $f$ may not be an internal ray, in what follows, we define a topological polynomial $F$ such that it preserves the internal rays and agrees with the behavior of $f$ on the closure of the basin of infinity.

Let us consider the restriction $f:\partial U\to \partial f(U)$. Then $g_U:=\phi_{f(U)}\circ f\circ\phi_U^{-1}:\partial \mb{D}\to\partial\mb{D}$ is a self-covering of the circle. We continuously extend $g_U$ into the whole disk $\mb{D}$ as 
$$g_U: re^{2\pi i\theta}\mapsto r^{\delta_U}\cdot g_U(e^{2\pi i\theta}),$$
where $\delta_U:=\tu{deg}(f|_U)$. Finally, the map $F$ associated to the parametrization $\Phi_f$ is defined by 
\begin{equation}\label{equ_def}
F:=
\left\{
\begin{array}{ll}
\phi_{f(U)}^{-1}\circ g_U\circ\phi_U &~~~\tu{on each bounded Fatou domain }U;\\
f  &~~~\tu{otherwise}.
\end{array}
\right.
\end{equation}
Since radial rays are preserved by $g_U$, the map $F$ sends internal rays of $U$ onto that of $F(U)$. Recall that $\Omega_f$ is the basin of infinity.
Based on Lemma \ref{lem:local_connectivity},
we have the following. 
\begin{lem}
	The above map $F:\mb{C}\to\mb{C}$  is a branched covering with $\tu{deg}(F)=\tu{deg}(f)$ and $F|_{\ol{\Omega}_f}=f|_{\ol{\Omega}_f}$. Thus $F$ is a topological polynomial.
\end{lem}
\begin{proof}
 	We first consider the continuity of  $F$. 
 	It suffices to check that $F$ is continuous at each $z$ in the Julia set $J_f$. 
	For each $\epsilon>0$, by the continuity of $f$, we have 
	\begin{equation}\label{eq:continuity}
	f(B(z,\delta_0))\subseteq B(f(z),\epsilon)
	\end{equation} for some $\delta_0>0$. 
	Let $W$ be the union of all bounded Fatou domains whose boundaries contain $z$. The number of such Fatou domains are finite. It is then clear that $F$ is continuous on $W$. Hence there exists $\delta_1>0$ such that $$F(B(z,\delta_1)\cap W)\subseteq B(f(z),\epsilon).$$
	By Lemma \ref{lem:local_connectivity}, one can also choose sufficiently small $\delta_2>0$ such that each bounded Fatou domain, which intersects $B(z,\delta_2)$ but is disjoint from $W$, is contained in $B(z,\delta_0)$. Then \eqref{eq:continuity} holds for $F$ and $B(z,\delta)$ with $\delta:=\tu{min}\{\delta_0,\delta_1,\delta_2\}$. 
	
	It remains to show that $F$ is a branched covering, that is, for each pair $z$ and $w:=F(z)$, there exist orientation preserving homeomorphisms $\zeta:N_z\to \mb{D}$ and $\xi: N_{w}\to\mb{D}$ on some neighborhoods $N_z$ and $N_{w}$ of $z$ and $w$, respectively, such that $\xi\circ F\circ \zeta^{-1}(x)=x^k.$ Let $\tu{Crit}(f)$ be the set of ciritical points of $f$. Write $$\tu{Crit}(F):=(\tu{Crit}(f)\cap J_f)\cup\{c_U: U \tu{ is a critical Fatou domain}\},$$
	and $\tu{CV}(F):=F(\tu{Crit}(F))$.
	Then $F$ is locally one-to-one at each $z\in \mb{C}\setminus\tu{Crit}(F)$. By the Domain Invariance Theorem, the map $$F:\mb{C}\setminus F^{-1}(\tu{CV}(F))\to \mb{C}\setminus \tu{CV}(F)$$ is actually a covering. 
	
	For each $z\in\tu{Crit}(F)$, 
	let $k$ be the degree of $F:N_z\setminus\{z\}\to N_{w}\setminus\{w\}$, where $N_{w}\setminus\{w\}$ is a punctured disk disjoint from $\tu{CV}(F)\cup \tu{CV}(f)$ and $N_z$ is the component of $F^{-1}(N_{w})$ containing $z$. Clearly 
	\begin{equation}\label{equ_def}\notag
	k:=\tu{deg}\left(F|_{N_{z}\setminus\{z\}}\right)=
	\left\{
	\begin{array}{ll}
	\tu{deg}(f|_{N_z}) &~~\tu{ if }z\in J_f;\\
	\tu{deg}(f|_U)  &~~\tu{ if }z=c_U.
	\end{array}
	\right.
	\end{equation}
	Let $\xi:N_{w}\to \mb{D}\ \ w\mapsto 0$ be a homeomorphism. Then the lift $\zeta$ of $\xi$ under the coverings $F|_{N_{z}\setminus\{z\}}$ and $x\mapsto x^k$ on $\mb{D}\setminus\{0\}$ can be extended to a homeomorphism from $N_z$ to $\mb{D}$. Hence $\xi\circ F\circ\zeta^{-1}(x)=x^k$ for all $x\in\mb{D}$. Similarly, one can check that at each $z\in F^{-1}(CV(F))\setminus \tu{Crit}(F)$ the above maps $\zeta$ and $\xi$ exist. Therefore, $F$ is a branched covering of degree $\tu{deg}(f)$,  whose critical set is $\tu{Crit}(F)$. The proof of the lemma is complete.
\end{proof}


\subsection{Regulated arcs}
Let $X$ be a topological space. Recall that $X$ is called \emph{arcwise connected} if each pair of distinct points in $X$ can be connected by an arc in $X$; and is called \emph{locally arcwise connected} (resp. \emph{locally connected}) if every point $z\in X$ has arbitrarily small arcwise connected (resp. connected) neighborhoods.
\begin{lem}[{\cite[Lemma 17.17 and 17.18]{Mi06}}]\label{lem:arcwise_collected}
	Every compact and locally connected metric space is locally arcwise connected.
\end{lem}

\begin{lem}\label{lem:Fatou domain}
   Suppose that the filled Julia set $K_f$ of $f$ is locally connected (or equivalently $J_f$ is locally connected). Then 
	\begin{enumerate}
		\item $K_f$ is arcwise connected and
		\item every component of $K_f\setminus\{p\}$ is arcwise connected for each $p\in J_f$.
	\end{enumerate}
\end{lem}
\begin{proof}
	We claim that a connected and locally arcwise connected set $A$ in $\mb{C}$ is arcwise connected. Indeed, after fixing a point $p\in A$, we set the non-empty set $Y\subseteq A$ as
	$$Y=\{p\}\cup\{z\in A: \tu{ there is an arc in $A$ joining $p$ to $z$}\}.$$
	Then both $Y$ and $A\setminus Y$ are open subsets of $A$. Since $A$ is connected, it holds that $A=Y$. Thus $A$ is arcwise connected.
	
	The statement of (1) follows directly from the claim and Lemma \ref{lem:arcwise_collected}. Now for any arcwise connected neighborhood $N\subseteq K_f$ of a point $z\in C$, if $p\notin N$, then $N\subseteq C$. Hence $C$ is locally arcwise connected. By the claim again, the statement of (2) follows.
	
\end{proof}

\begin{defi}[Regulated and clean arcs]\label{df:regulated}
	Let $\gamma$ be an arc in $K_f$ with both endpoints in the Julia set $J_f$. The arc $\gamma$ is  called \emph{regulated}, if for each bounded Fatou domain $U$ whose closure intersects $\gamma$, the intersection $\gamma\cap \ol{U}$ is 
	$$\emph{either a singleton or the union of two closed internal rays}.$$
	The arc $\gamma$ is called \emph{clean} if $\gamma$ intersects the closure of each bounded Fatou domain in at most one point. Clearly clean arcs are always regulated.
\end{defi}

The following lemma establishes strong rigidity of regulated arcs: the existence and uniqueness.
\begin{lem}\label{lem:regulated_arcs}
	Let $F$ be the topological polynomial associated to $f$. Then for all $x\neq y\in J_f$, we have:
	\begin{enumerate}
		\item There exists a unique regulated arc $\gamma$ connecting $x$ and $y$. We write $\gamma$ as $[x,y]$, the open arc $\gamma\setminus\{x,y\}$ as $]x,y[$, and $\gamma\setminus\{x\}$ as $]x,y]$.
		\item The image $F([x,y])$ contains $[f(x),f(y)]$ and has no loops (thus is a tree).
		\item If $]x,y[$ is disjoint from $\tu{Crit}(F)$, then $F:[x,y]\to [f(x),f(y)]$ is a homeomorphism.
	\end{enumerate}
\end{lem}
Before proving the lemma, we introduce some terminologies. Let $\gamma(t):[0,1]\to \mb{C}$ be an arc and $E$ be a subset of $\mb{C}$ such that $\gamma\cap E\neq\emptyset$. We set the \emph{first-in} (resp. \emph{last-out}) place of $\gamma$ meeting $E$ as $\gamma(t_1)$ (resp. $\gamma(t_2)$), where
$$t_1:=\tu{inf}\,\{t\geq 0: \gamma(t)\in\gamma\cap E\} \tu{ and } t_2:=\tu{sup}\,\{t\geq 0: \gamma(t)\in \gamma\cap E\}.$$
The case $\gamma(t_1)=\gamma(t_2)$ happens if and only if $\gamma\cap E$ is a singleton.
\begin{proof}
	(1) By Lemma \ref{lem:local_connectivity}, one may enumerate all bounded Fatou domains of $f$ as $U_n$ for $n\geq 1$, such that $$\tu{diam}\,U_n\geq \tu{diam}\,U_{n+1}.$$
	Since $K_f$ is arcwise connected, let $\gamma_0$ be an arc in $K_f$ connecting $x$ and $y$. We will inductively construct a sequence of arcs $\gamma_n$ such that their limit $\gamma$ is as required.
	
	For $n\geq 1$, if $\gamma_{n-1}\cap \partial U_n$ contains at most one point, then we set $\gamma_n=\gamma_{n-1}$. Otherwise, let $x_n$ and $y_n$ be the first-in and last-out places of $\gamma_{n-1}$ meeting $\partial U_n$, respectively. We first remove the open segment of $\gamma_{n-1}$ bounded by $x_n$ and $y_n$, then replace it by the union of the two internal rays of $U_n$ landing at $x_n$ and $y_n$. The new arc is denoted by $\gamma_n$.
	
	By induction, we obtain a sequence of arcs $\gamma_n$ for $n\geq 0$. Note that $\gamma_{n}$ differs from $\gamma_{n-1}$ only possibly in the Fatou domain $U_{n}$. As a consequence of the shrinking of the sizes of $U_n$ as $n$ tends to $\infty$, the sequence $\gamma_n$ converges uniformly to a curve $\gamma$. From the construction, $\gamma$ is a regulated arc linking $x$ and $y$. 
	
	For the uniqueness, we prove it by contradiction and assume that $\gamma'$ is another regulated arc connecting $x$ and $y$. Then there is a bounded component $W$ of the set $\mb{C}\setminus(\gamma\cup\gamma')$. Clearly $W$ is disjoint from $\Omega_f$. So $W\cap \ol{\Omega}_f=\emptyset$. Then $W$ is contained in a bounded Fatou domain, say $U$. We have
\begin{equation}\label{eq:regulated}
\partial W=\partial W\cap \ol{U}\subseteq (\gamma\cup \gamma')\cap \ol{U}=(\gamma\cap \ol{U})\cup (\gamma'\cap \ol{U}).
\end{equation}
	By definition, the term on the right hand side of \eqref{eq:regulated} contains at most four closed internal rays of $U$, while $\partial W$ is a Jordan curve. This is impossible.
	
	$(2)$ The image $F([x,y])$ is connected. It intersects the Fatou set in several internal rays. Moreover, $F([x,y])\cap \partial U$ is a finite set for each bounded Fatou domain $U$. 
	Let $\gamma'$ be an arc in $F([x,y])$ linking $f(x)$ and $f(y)$. In fact, $\gamma'$ is the regulated arc $\gamma=[F(x),F(y)]$. Otherwise, $\gamma\cup \gamma'$ bounds a Jordan domain $W$.  Then an argument similar to $\eqref{eq:regulated}$ on $\partial W$ gives a contradiction.
	
	$(3)$ It holds from the fact that $F$ is locally injective on $[x,y]$ and that the image $F([x,y])$ cannot possess loops.	
\end{proof}

Recall that a point in the Julia set is called \emph{branched}, if it is the landing point of at least three rays. Preimages of branched points are branched as well. Thus they (if exist) form a dense subset of the Julia set.
\begin{lem}\label{lem:branched_dense}
	Suppose that the locally connected Julia set $J_f$ is not a segment. Then the set of branched points is dense in each clean arc.
\end{lem}
\begin{proof}
Since subarcs of clean arcs are still clean, 
it suffices to show that a clean arc $\gamma=[x,y]$ contains at least one branched point. 


Let $W\subseteq J_f$ be an arcwise connected neighborhood of a point $w\in]x,y[$ such that $x,y\notin W$. We claim that the set $W\setminus \gamma$ contains at least one point in the Julia set. Indeed, if $f$ has no bounded Fatou domains, as $J_f$ is not a segment, the branch points are dense in $W$. Thus such a point $z$ exists. Otherwise, 
the existence of clean arcs implies that there are infinitely many Fatou domains. The union of the boundaries of them is dense in $J_f$. So such a point $z$ exists.   
	
Let $\gamma_{wz}$ be an arc (maybe not regulated) in $W$ joining $w$ and $z$. Let $p$ be the first-in place of $\gamma_{wz}$ meeting $]x,y[$. Let $\gamma_{zp}$ be the subarc of $\gamma_{zw}$ joining $z$ and $p$. 

Starting with $\gamma_{zp}$, one can obtain a regulated arc $[z,p]$, according to the proof of Lemma \ref{lem:regulated_arcs} (1), such that $$[z,p]\cap J_f\subseteq \gamma_{zp}.$$
Then $[z,p]\cap[x,y]=\{p\}$. The three arcs $[x,p],[y,p]$ and $[z,p]$ form a ``Y" shape, precisely,
	$$[x,p]\cap[y,p]=[z,p]\cap [x,p]=[z,p]\cap[y,p]=\{p\}.$$
	
It suffices to show that $x,y$ and $z$ belong to distinct components of $K_f\setminus\{p\}$ by \cite[Theorem 6.6]{Mc95}. For otherwise, two of them, say $x$ and $z$, can be linked by an arc $\gamma_{xz}$ in a component of $K_f\setminus\{p\}$ by Lemma \ref{lem:Fatou domain}. Again one can derive the regulated arc $[x,z]$ from $\gamma_{xz}$. It holds that $[x,z]\cap J_f\subseteq \gamma_{xz}\cap J_f$. In particular, $p\notin[x,z]$ as $p\notin \gamma_{xz}$. Thus the union $[x,p]\cup[p,z]\cup [x,z]$ possesses a loop, which is a contradiction. The proof is complete.
\end{proof}

\section{Critical portraits}\label{critical_portraits}
In this section, we first introduce the critical portraits of a polynomial and then study the properties of them. A critical portrait induces a partition of the dynamical plane $\mb{C}$. On each piece of the partition, the behavior of the associated topological polynomial $F$ (defined in \eqref{equ_def}) can be well understood (Proposition \ref{prop:homeo_on_pieces} and \ref{prop:partition_c}). Throughout this section, the Julia set $J_f$ is assumed to be locally connected.

The following notions of supporting rays and angles are needed to define the critical portraits. 
\begin{defi}
	Consider a bounded Fatou domain $U$ and a point $z\in\partial U$. The union of rays landing at $z$ separates $\mb{C}\setminus\{z\}$ into several parts. One of them that contains $U$ is assumed to be bounded by $R(\theta_1)$ and $R(\theta_2)$. Then $R(\theta_1)$ and $R(\theta_2)$ (resp. $\theta_1$ and $\theta_2$) are called the rays (resp. angles) \emph{supporting} at $(U,z)$. The trivial case that $\theta_1=\theta_2$  happens if and only if exactly one ray lands at $z$. 	
\end{defi}

In order to create a critical portrait, we assign some subsets $\Theta(U)$ and $\Theta(c)$ of $\mathbb{T}$ to each critical Fatou component $U$ and each critical point $c\in J_f$ of $f$ as follows.
\begin{enumerate}
	\item[(A1)] For a critical point $c\in J_f$, take a ray landing at $f(c)$. By pulling this ray back via $f$, we obtain $\delta_c:=\tu{deg}(f|_c)$ rays at $c$. We define $\Theta(c)$ as the collection of the angles of these $\delta_c$ rays.
	
	
	\item[(A2)] For a strictly pre-periodic critical Fatou domain $U$, choose a ray $R$ that supports $f(U)$ at a point $w\in\partial f(U)$ (maybe a critical value). Since $f:\partial U\to \partial f(U)$ is a covering of degree $\delta_U:=\tu{deg}(f|_U)$, there are $\delta_U$ preimages $z_1,\cdots,z_{\delta_U}$ of $w$ in $\partial U$. For each $z_k$, choose a ray $R_k$ supporting at $(U,z_k)$ such that $f(R_k)=R$ \footnote{There are two choices of such $R_k$ if and only if $z_k$ is a critical point and just one ray lands at $w$; see $\Theta(c_2)$ and  $\Theta(U_3)$ in Figure \ref{fig:partition}.}.
	For convenience, we let all $R_k$ support $U$ at the same direction.
	We set $\Theta(U)$ as the collection of the angles of $R_1,\cdots,R_{\delta_U}$.
	\item[(A3)] For a period $p$ cycle of Fatou domains $U_0,U_1=f(U_0),\cdots, U_0=f(U_{p-1})$, we first choose a point $z_0\in\partial U_0$ fixed by $f^p$, and then pick a ray $R_0$ supporting at $(U_0,z_0)$. In the cycle 
	$$(U_0, z_0, R_0)\mapsto\cdots\mapsto (U_{p-1},z_{p-1},R_{p-1})\mapsto (U_p,z_p,R_p)(=(U_0,z_0,R_0)),$$
the rays $R_0,\cdots, R_{p-1}$ support at $(U_0,z_0),\cdots,(U_{p-1},z_{p-1})$, respectively, at the same direction. For each critical $U_k$ in the cycle, the way of setting $\Theta(U_k)$ is analogous to the procedure in (A2). By pulling back $(z_{k+1},R_{k+1})$ via $f|_{\partial{U_k}}$, we obtain $\delta_{U_k}$ preimages of $z_{k+1}$ in $\partial U_k$ and their corresponding supporting rays. We use $\Theta(U_k)$ to denote the $\delta_{U_k}$ angles of these supporting rays.
\end{enumerate}

With all the settings above, we let $\mc{C}_f$ be the family of these $\Theta(c)$ and $\Theta(U),$ where $c$ and $U$ are taken over all critical points and critical Fatou domains respectively. 


Two elements $\Theta'$ and $\Theta''$ in $\mc{C}_f$ are said to be \emph{equivalent}, if either $\Theta'\cap\Theta''\neq \es$ or there are $\Theta_1,\cdots,\Theta_k$ in $\mc{C}_f$ linking $\Theta'$ and $\Theta''$ in the sense that
$$\Theta'\cap\Theta_1\neq \es,\Theta_1\cap\Theta_2\neq\es,\cdots,\Theta_k\cap\Theta''\neq \es.$$
Each pair of distinct elements in $\mc{C}_f$ are clearly disjoint when $f$ has no bounded Fatou domains.
\begin{defi}[Critical portrait]\label{def:critical_portrait}
For each $\Theta\in\mc{C}_f$, we denote by $\wt{\Theta}$ the union of all elements in $\mc{C}_f$ that are equivalent to $\Theta$. Then the collection $\wt{\mc{C}}_f=\{\wt{\Theta}_1,\cdots,\wt{\Theta}_m\}$ is called a \emph{critical portrait} (of $f$). See Figure \ref{fig:partition}. 	
\end{defi}
As a straightforward consequence of this construction, critical portraits obey the following rules:
\begin{enumerate}
	\item[(R1)] For each $\wt{\Theta}_k$, its image under $m_d:\theta\mapsto d\theta\tu{ mod }\mb{Z}$ contains exactly one angle.
	\item[(R2)] $\sum_{1\leq k\leq m}(\#\wt{\Theta}_k-1)=d-1$, due to Hurwitz's formula.
	\item[(R3)] $\wt{\Theta}_1,\cdots,\wt{\Theta}_m$ are pairwise \emph{unlinked}, that is, the sets $\wt{\Theta}_k$ and $\wt{\Theta}_{k'}$ are contained in two disjoint intervals of $\mb{T}$ for each pair $k\neq k'$.
\end{enumerate}

\begin{rmk}\label{rmk:critical_portrait}
	The notion of critical portraits defined here is slightly broader than that in \cite{Po09}, as the set $\Theta(U)$ in (A2) does not need to satisfy the iterated condition in \cite{Po09}.
\end{rmk}

\subsection{Partitions of $\ol{\mb{D}},\mb{T}$ and $\mb{C}$}
\begin{figure}[h]
	\begin{tikzpicture}
	\node at (0,0) {\includegraphics[width=15cm]{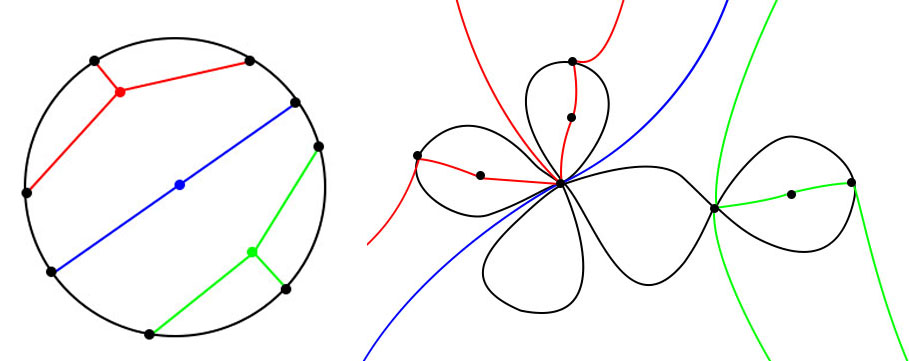}};
	\node at (0.25,-0.25) {$U_1$};
	\node at (2.25,1.0) {$U_2$};
	\node at (6.,-0.5) {$U_3$};
	\node at (2.15,-0.20) {$c_1$};
	\node at (4.5,-0.75) {$c_2$};
	\node at (-1, -1){$R(\theta_1)$};
	\node at (-0.5,2.5){$R(\theta_2)$};
	\node at (2.3,2.5){$R(\theta_3)$};
	\node at (-0.5,-2.5){$R(\theta_4)$};
	\node at (3.75,2.5){$R(\theta_5)$};
	\node at (5.5,2.5){$R(\theta_6)$};
	\node at (4.75,-2.5){$R(\theta_7)$};
	\node at (7,-2.5){$R(\theta_8)$};
	\end{tikzpicture}
	\caption{The critical Fatou domains are $U_i, i=1,2,3,$ and the critical points are $c_1$ and $c_2$ at the boundaries. As shown, $\Theta(U_1)=\{\theta_1,\theta_2\}, \Theta(U_2)=\{\theta_2,\theta_3\}, \Theta(U_3)=\{\theta_6,\theta_8\}$, $\Theta(c_1)=\{\theta_4,\theta_5\}$ and $\Theta(c_2)=\{\theta_6,\theta_7\}$. By definition, the critical portrait $\wt{C}_f=\{\wt{\Theta}_1,\wt{\Theta}_2,\wt{\Theta}_3\}$ with $\wt{\Theta}_1=\Theta(U_1)\cup\Theta(U_2)$, $\wt{\Theta}_2=\Theta(c_1)$ and $\wt{\Theta}_3=\Theta(c_2)\cup\Theta(U_3)$. Associated to $\mc{C}_f$, the partitions $\mc{I}_f$ and $\mc{D}_f$ are shown in the left picture; the right picture shows $\mc{P}_f=\{P_1,\cdots,P_6\}$ of $\mb{C}$. }
	\label{fig:partition}
\end{figure}

A critical portrait naturally induces partitions: $\mc{D}_f,\mc{I}_f$, and $\mc{P}_f$ of the closed unit disk $\ol{\mb{D}}$, the unit circle $\mb{T}$, and the plane $\mb{C}$, respectively; see Figure \ref{fig:partition}. We describe the precise constructions as follows.

In the boundary $\partial \mb{D}$, for each $\wt{\Theta}_k$, we first mark all the points $e^{2\pi i\theta}$ for $\theta\in \wt{\Theta}_k$; then draw $\#\wt{\Theta}_k$ straight line segments, each of which starts at a marked point and ends at the center of gravity of the $\#\wt{\Theta}_k$ marked points; the union of these closed segments is denoted by $Y_k$. Then by rules (R2) and (R3), one has
\begin{itemize}
	\item $Y_k$ and $Y_{k'}$ are disjoint for $k\neq k'$;
	\item $Y_1,\cdots,Y_m$ cut $\ol{\mb{D}}$ into $d$ pieces $D_1,\cdots,D_d$.
\end{itemize} 
The collection $\mc{D}_f=\{D_1, \cdots, D_d\}$ is the partition of $\ol{\mb{D}}$ induced by $\wt{\mc{C}}_f$.

In the unit circle $\mb{T}$, associated to each $D_k$ above, there is an open subset $I_k$ given by 
$$I_k=\{\theta\in\mb{T}: e^{2\pi i\theta}\in D_k\cap\partial\mb{D}\}.$$
The collection $\mc{I}_f=\{I_1,\cdots,I_d\}$ of $I_k$ is a partition of $\mb{T}$. And again, from rules (R1) and (R2), we have
\begin{itemize}
	\item the total length of each $I_k$ is $1/d$;
	\item the map $m_d:I_k\to \mb{T}\setminus m_d(\partial I_k)$ is bijective.
\end{itemize} 

To get the corresponding partition of $\mb{C}$, we need some further notations as follows:
\begin{enumerate}
	\item[-] $\mc{R}(c)$: the union of the critical point $c$ and all $R(\theta)$ for $\theta\in \Theta(c)$;
	\item[-] $\mc{R}(U)$: the union of all $R(\theta)$, their distinct landing points $z_\theta\in\partial U$ and the internal rays of $U$ landing at $z_\theta$ for $\theta\in \Theta(U)$;
	\item[-] $\wt{\mc{R}}_k$: write $\wt{\Theta}_k$ as 
	\begin{equation}\label{eq:form}
	\wt{\Theta}_k=\Theta(c_1)\cup\cdots\Theta(c_l)\cup\Theta(U_1)\cup\cdots\cup\Theta(U_{l'}),
	\end{equation}
	and then $\wt{\mc{R}}_k$ is the union of $$\mc{R}(c_1),\cdots,\mc{R}(c_l),\mc{R}(U_1),\cdots,\mc{R}(U_{l'}).$$
	\item[-] $\wt{\mc{R}}_f:=\{\wt{\mc{R}}_1,\cdots,\wt{\mc{R}}_m\}$ for a critical portrait $\wt{\mc{C}}_f=\{\wt{\Theta}_1,\cdots,\wt{\Theta}_m\}$.
\end{enumerate}

\begin{lem}\label{lem:lines}
The elements in $\wt{\mc{R}}_f$ satisfy the followings:
	\begin{enumerate}
		\item $\wt{\mc{R}}_k\cap K_f$ is connected and is a tree.
		\item When $\wt{\mc{R}}_i\cap \wt{\mc{R}}_j\neq \emptyset$ with $i\neq j$, the intersection is a singleton within the Julia set.
		\item $\wt{\mc{R}}_i$ never crosses $\wt{\mc{R}}_j$, i.e., $\wt{\mc{R}}_i$ is contained in the closure of a component of $\mb{C}\setminus\wt{\mc{R}}_j$.
		\item Let $\theta:=m_d(\wt{\Theta}_k)$ and $z_\theta$ be the landing point of $R(\theta)$. Then the image $F(\wt{\mc{R}}_k)$ is one of the following types:
		\begin{enumerate}
			\item[type I:] the closed ray $\ol{R(\theta)}$.
			\item[type II:] the closed ray $\ol{R(\theta)}$ together with an internal ray landing at $z_\theta$.
			\item[type III:] the closed ray $\ol{R(\theta)}$ together with two internal rays (they are from distinct Fatou domains but land together at $z_\theta$). 
		\end{enumerate}
	\end{enumerate}
\end{lem}
\begin{proof}
	$(1)$ Clearly $\mc{R}_k\cap K_f$ is connected. An argument similar to the proof of Lemma \ref{lem:regulated_arcs} (3) yields that there is no loop in $\mc{R}_k\cap K_f$, i.e., it is a tree.
	
	$(2)$ By definition, each bounded Fatou domain can not have non-empty intersection with both $\wt{R}_i$ and $\wt{R}_j$. Thus $\wt{R}_i\cap \wt{R}_j$ belongs to the Julia set. Then it is a finite set. If it contains more than one point, then there are loops in $\wt{R}_i\cup\wt{R}_j$, which is impossible.
	
	
	$(3)$ The statement follows from the fact that $\wt{\Theta}_i$ and $\wt{\Theta}_j$ are unlinked; see (R3).
	
	$(4)$ We write $\wt{\Theta}_k$ in the form of \eqref{eq:form}.
	From the construction, $R(\theta)$ is the only external ray in $F(\wt{\mc{R}}_k)$. Moreover, 	
	$F(\wt{\mc{R}}_k)=\ol{R(\theta)}$ (is of type I) if and only if $l'=0$ (then $l=1$). 
	
	If $l'\geq 1$, then $R(\theta)$ supports the Fatou domains $F(U_1),\cdots,F(U_{l'})$ at $z_\theta$. Since each ray supports at most two Fatou domains, the number of $F(U_1),\cdots,F(U_{l'})$ is either one or two, which depends on the type (II or III) of $F(\wt{\mc{R}})$.
\end{proof}

Let us consider a pair of associated pieces $(D,I)$ from $(\mc{D}_f,\mc{I}_f)$. The set $\partial D\setminus\partial{\mb{D}}$ is a disjoint union of open segments, say $l_{\theta_1\theta_1'},\cdots,l_{\theta_n\theta_n'}$, in $\mb{D}$. The endpoints of each $l_{\theta_j\theta_j'}$ are denoted by $e^{2\pi i\theta_j}$ and $e^{2\pi i\theta'_j}$. Let $\mathcal{L}_{\theta_j\theta_j'}\subseteq\mb{C}$ be the arc consisting of
\begin{enumerate}
	\item the rays $R(\theta_j)$ and $R(\theta_j')$ together with their landing points $z_{\theta_j}$ and $z_{\theta_j'}$;
	\item the regulated arc $[z_{\theta_j},z_{\theta_j'}]$ (it is contained in $\wt{\mc{R}}_j\cap K_f$).
\end{enumerate}

Each $\mathcal{L}_{\theta_j\theta_j'}$ cuts the plane $\mb{C}$ into two components. Let $\mathcal{L}^+_{\theta_j\theta_j'}$ be the component containing $R(\theta)$ for $\theta\in I$. Then the intersection 
$$P:=\mc{L}^+_{\theta_1\theta_1'}\cap\cdots\cap \mc{L}^+_{\theta_n\theta_n'},$$
which is open and might be disconnected, is the piece associated to $(D,I)$. See Figure \ref{fig:partition}. 

Let $\mc{P}_f$ be the collection of these $d$ pieces associated to $(\mc{D}_f,\mc{I}_f)$. Then $\mc{P}_f$ is the partition of the dynamical plane $\mb{C}$ induced by $\wt{\mc{C}}_f$.

We emphasize that elements in $\mc{P}_f$ are one to one corresponding to elements in $\mc{D}_f$ and those in $\mc{I}_f$. The families $\mc{D}_f$ and $\mc{I}_f$ are used as models to look at the pieces in $\mc{P}_f$.

The next lemma describes some properties of $\mc{P}_f$. The proof is omitted, since it follows directly from Lemma \ref{lem:lines} and the construction of the partitions.
\begin{lem}
	Consider $\mc{P}_f=\{P_1,\cdots,P_d\}$ induced by a critical portrait $\wt{\mc{C}}_f$. Then
	\begin{enumerate}
		\item each $P_i$ is open and the components of $P_i$ are unbounded;
		\item $P_i\cap\Omega_f=\cup_{\theta\in I_i}R(\theta)$ with $I_i\in\mc{I}_f$ corresponding to $P_i$; 
		\item $P_i\cap P_j=\es$ whenever $i\neq j$;
		\item $\mb{C}\setminus(\wt{\mc{R}}_1\cup\cdots\cup\wt{\mc{R}}_m)=P_1\cup\cdots\cup P_d$.
	\end{enumerate}
\end{lem}

\subsection{The behavior of $F$ on pieces}
\begin{prop}\label{prop:homeo_on_pieces}
	Let $P$ be a piece in $\mc{P}_f$. Then $F:P\to \mb{C}\setminus F(\partial P)$ is bijective.
\end{prop}
\begin{proof}
	The proof is based on the Topological Argument Principle. It states that the number of solutions, counted with multiplicity, to the equation $$F(z)=z_0\tu{ with }z\in W,$$  is the winding number of $F(\partial W)$ around $z_0$, where $W$ is an open and bounded set in $\mb{C}$ with $\partial W$ consisting of finitely many arcs. We denote the winding number by $\tu{wind}(F(\partial W),z_0)$.
	
	Recall that $\Psi_f:\Omega_f\to\mb{C}\setminus\ol{\mb{D}}$ is the B\"{o}ttcher map. For $t>1$, let $G_t$ be the bounded domain surrounded by the Jordan curve $\Psi_f^{-1}(\{\zeta\in\mb{C}: |\zeta|=t\})$. 
	
	For each $t>1$, we consider the open set $P_t:=P\cap G_t$, whose boundary consists of edges of two types:
	\begin{enumerate}
		\item $\Psi_f^{-1}(\{te^{2\pi i\theta}:\theta\in I\})(=P\cap \partial G_t)$, which is mapped bijectively onto $$\partial G_{t^d}\setminus \Psi_f^{-1}(\{t^de^{2\pi im_d(\theta)}: \theta\in \partial I\});$$
		\item $G_t\cap \mc{L}_{\theta_1\theta_1'},\cdots, G_t\cap \mc{L}_{\theta_n\theta_{n'}}$, whose images under $F$ are formed by closed internal rays and sub-arcs of external rays by Lemma \ref{lem:lines}.
	\end{enumerate}
	Then we conclude that
	\begin{equation}\notag
	\tu{wind}(F(\partial P_t),z_0)=
	\left\{                       
	\begin{array}{ll}              
	1 &~~~~\text{if}~z_0\in G_{t^d}\setminus F(\partial P), \\           
	0 &~~~~\text{if}~z_0\in\mb{C}\setminus\ol{G_{t^d}}.  
	\end{array}                  
	\right.                      
	\end{equation}
Thus $F:P_t\to G_{t^d}\setminus F(\partial P)$ is one-to-one. 	
	By the arbitrariness of $t$, the proof is complete.	
\end{proof}
\begin{prop}\label{prop:partition_c}
	Let $P$ be a piece in $\mc{P}_f$. For all $x\neq y$ in $P\cap J_f$, we have
	\begin{enumerate}
		\item $[x,y]\subseteq \ol{P}$.
		\item $F:[x,y]\to [f(x),f(y)]$ is bijective.
	\end{enumerate}
\end{prop}
\begin{proof}
	$(1)$ It suffices to prove that $[x,y]$ is contained in $\mc{L}^+_{\theta_j\theta_j'}$ for each $1\leq j\leq n$. This holds clearly when $\#[x,y]\cap\mc{L}_{\theta_j\theta_j'}\leq 1$. Otherwise,
	there exist distinct $z_1,z_2\in [x,y]$ such that $[x,y]$ meets $\mc{L}_{\theta_j\theta_j'}$ at the first-in place $z_1$ and last-out place $z_2$. Then $z_1$ and $z_2$ are either Fatou centers or points in the Julia set. One may break $[x,y]$ into three segments $$[x,y]=\gamma_{xz_1}\cup\gamma_{z_1z_2}\cup\gamma_{z_2y}.$$ The absence of loops in $\mc{L}_{\theta_j\theta_j'}\cup[x,y]$ implies that $\gamma_{z_1z_2}\subseteq \mc{L}_{\theta_j\theta_j'}$. By Lemma \ref{lem:lines} (3), the arcs $\gamma_{xz_1}$ and $\gamma_{z_2y}$ are contained in the same component of $\mb{C}\setminus\mc{L}_{\theta_j\theta_j'}$, which must be $\mc{L}^+_{\theta_j\theta_j'}$. 
	
	$(2)$ If $F:[x,y]\to [f(x),f(y)]$ is not bijective, we assume $z_1\neq z_2\in[x,y]$ such that $w=F(z_1)=F(z_2)$. 	
	If $w$ is a point in the Fatou set but not a Fatou center, then each $z_k$ belongs to the interior of an internal ray $R_k$ of some Fatou domain $U_k$. Then $$R_1\neq R_2\subseteq [x,y]\cap \partial P\tu{ and }F(R_1)=F(R_2).$$	
	
	If $U_1=U_2$, then $R_1,R_2\subseteq\mc{R}(U_1)$. The set $[x,y]\setminus U_1$ has two  components. Both of them are non-trivial regulated arcs contained in the same component of $\mb{C}\setminus \mc{R}(U_1)$.
	This cannot happen, because the rays in $\mc{R}(U_1)$ supporting at $(U_1, z_1)$ and $(U_1,z_2)$ have the same direction. 
	
	If $U_1\neq U_2$, pick a neighborhood $N_k$ of $z_k$. Let $N_k'=N_k\cap P$. Since $F|_{P}$ preserves the orientation, $F(N_1')$ and $F(N_2')$ have non-empty intersection. Thus $F|_P$ can not be injective, which contradicts Proposition \ref{prop:homeo_on_pieces}.
	 
	We conclude that $w$ is either a Fatou center or a point in $F(\partial P)\cap J_f$. Then $F([x,y])$ has only finitely many self-intersection points. The absence of loops in $F([x,y])$ gives that such a self-intersection point $w$ does not exist. The proof is complete. 
\end{proof}
\begin{defi}\label{itinerary}
	Let $f$ be a polynomial of degree $d\geq 2$ with locally connected Julia set. Let $\wt{\mc{C}}_f$ be a critical portrait of $f$. Let $\mc{I}_f:=\{I_1,\cdots,I_d\}$ and $\mc{P}_f:=\{P_1,\cdots,P_d\}$ be the partitions of $\mb{T}$ and $\mb{C}$ associated to $\wt{\mc{C}}_f$. For an angle $\theta\in\mb{T}$, whose forward orbit under $m_d$ avoids $\partial I_1\cup\cdots\cup\partial I_d$, its \emph{itinerary} is defined as
	$$\tu{itin}(\theta)=n_0n_1\cdots n_k\cdots \textup{ with }m_d^k(\theta)\in I_{n_k}.$$
	Similarly, the \emph{itinerary} of a point $z$ in the Julia set $J_f$, whose forward orbit is disjoint from $E_f:=(\partial P_1\cup\cdots\cup \partial P_d)\cap J_f$, is the sequence 	
	 $$\tu{itin}(z)=n_0n_1\cdots n_k\cdots\tu{ with } f^k(z)\in P_{n_k}.$$
	 If such a point $z$ is the landing point of a ray $R(\theta)$, then clearly $\tu{itin}(\theta)=\tu{itin}(z)$.
\end{defi}
Finally, we note that the itineraries of all but countably many elements in $\mb{T}$ and $J_f$ are well-defined. 

\section{Proof of the main theorem}\label{proof_main_thm}
With all the preparations in the previous sections, in this section we prove the main result Theorem \ref{thm:main_thm} of this paper. It states that the landing points of $R(\theta_1)$ and $R(\theta_2)$, for angles $\theta_1$ and $\theta_2$ with the same itinerary, must either coincide or lie in the boundary of a Fatou domain, which is eventually iterated onto a Siegel disk.


A clean arc $I$, defined in Definition \ref{df:regulated}, is called \emph{wandering} if $f^{i}(I)\cap f^{j}(I)=\es$ for all integers $i\neq j\geq 0$. Before proving Theorem \ref{thm:main_thm}, we need a few lemmas and propositions.
\begin{lem}\label{lem:non_wandering_arc}
	There is no wandering clean arc.
\end{lem}
\begin{proof}
	The lemma holds immediately when the Julia set $J_f$ is a segment. Otherwise,
	we assume that $I$ is a wandering clean arc. By replacing $I$ with some of its iterate when necessary, one may further assume that the forward orbit of $I$ is disjoint from all critical points of $f$. Then $f^n$ restricted on $I$ is injective for all $n\geq 1$.
	
	By Lemma \ref{lem:branched_dense}, branched points are dense in $I$. 
	All of them are wandering as $I$ is assumed to be wandering. From Corollary \ref{cor:finite_branched}, there exist two branched points in $I$, say $z_1$ and $z_2$, such that
	$$f^i(z_1)=f^j(z_2)\tu{ for some $i\geq 0$ and $j\geq 0$}.$$
	As $f^i|_{I}$ is injective, it holds that $i\neq j$.
	Thus $I$ is not wandering, contradicting our assumption. We conclude the lemma.
\end{proof}
\begin{prop}\label{prop:semi-buried}
	Let $[x,y]$ be a clean arc such that the itineraries of $x$ and $y$ are well-defined with respect to a critical portrait $\wt{\mc{C}}_f$. Then $\tu{itin}(x)\neq \tu{itin}(y)$.
\end{prop}
\begin{proof}
	We prove it by contradiction. Let $x_k=f^k(x)$ and $y_k=f^k(y)$ for all $k\geq 0.$
	We assume $\tu{itin}(x)=\tu{itin}(y)=s_0s_1\cdots s_k\cdots.$ Then by Proposition \ref{prop:partition_c}, for each $k\geq 1$, there exists a bijection
	$$f^k:[x_0,y_0]\to [x_k,y_k]\,\subseteq \ol{P_{s_k}}.$$
	Recall that $E_f:=(\partial P_1\cup \cdots\cup\partial P_d)\cap J_f$; see Definition \ref{itinerary}. According to Lemma \ref{lem:non_wandering_arc}, there is a point $\xi_0$ in the Julia set $J_f$ and two integers $l\geq 0$ and $n\geq 1$ such that
		\begin{enumerate}
		\item[$\bullet$] $\xi_0\in [x_{l+n},y_{l+n}]\cap[x_l,y_l]$;
		\item[$\bullet$] the forward orbit of $\xi_0$ is disjoint from $E_f$ as $\#E_f<\infty$. 
	\end{enumerate}
	By replacing $[x_l,y_l]$ with $[x_0,y_0]$, one may assume $l=0$. Let us consider the set
	$$H:=[x_0,y_0]\cup [x_{n},y_{n}]\cup\cdots\cup[x_{kn},y_{kn}]\cup\cdots.$$
 	Since the successive arcs $[x_{kn},y_{kn}]$ and $[x_{(k+1)n},y_{(k+1)n}]$ for $k\geq 0$, possess a common point $\xi_{kn}:=f^{kn}(\xi_0)\in P_{s_{kn}}$, it follows that
 	$$s_0=s_{n}=\cdots =s_{kn}=\cdots.$$
 	Thus $H\subseteq \ol{P_{s_0}}$. Moreover, as $H$ is arcwise connected and has no loops by the rigidity of regulated arcs, $H$ is a finite or infinite tree.  	
 	Since $f$ is one-to-one on $P_{s_0}$, the restricted map $f^{n}|_H$ is injective with the exception of finitely many points in $\partial P_{s_0}$. Again as $f^{n}(H)\subseteq H$ has no loops, the map
 	$g:=f^{n}:H\to H$ is injective. 
	
	We claim that $g$ has no fixed points. For otherwise, there exists some $z$ with $g(z)=z$. Then $z\in [x_0,y_0]$. There exists a small subarc $[z,z_\epsilon]$ of $[x_0,y_0]$ such that one of the three cases occurs: 
	$$[z,z_\epsilon]\subseteq [z,g(z_\epsilon)[, ~~~[z,z_\epsilon[~\supseteq [z,g(z_\epsilon)]\tu{ and }[z,z_\epsilon]\cap[z,g(z_\epsilon)]=\{z\}.$$ 
	We now obtain contradictions as follows. 
	\begin{enumerate}
		\item If $[z,z_\epsilon]\subseteq [z,g(z_\epsilon)[$, then the arcs 
		$$g^k\,(\,] z_\epsilon,g(z_\epsilon)[\,)=[g^{k+1}(z_\epsilon),g^k(z_\epsilon)]$$ for $k\geq 0$ are pairwise disjoint by the injective property of $g$ on $H$. Thus $]z_\epsilon,g(z_\epsilon)[$ is wandering under $g$. This contradicts Lemma \ref{lem:non_wandering_arc}.
		
		\item If $[z,z_\epsilon]\supseteq [z,g(z_\epsilon)]$, which means that $g$ is attracting on $[z,z_\epsilon]$ and $$]z,z_\epsilon]=\cup_{k\geq 0}[g^{k+1}(z_\epsilon),g^k(z_\epsilon)],$$ then $]g(z_\epsilon),z_\epsilon[$ is wandering under $g$. Again this is a contradiction. 
		
		\item If $[z,z_\epsilon]\cap[z,g(z_\epsilon)]=\{z\}$, then the absence of wandering clean arcs implies that, for some $k$, we have $$g^k\,(\,]z,z_\epsilon])\,\cap\, ]z,z_\epsilon]\neq \es.$$
		Since $H$ contains no loops, the arcs $g^k(\,[z,z_\epsilon])$ and $[z,z_\epsilon]$ must overlap on a subarc $[z,\wt{z}_\epsilon]$. By discussions on $[z,\wt{z}_\epsilon]$ and $g^k$ similar to Case (1) and (2), there is a contradiction.
	\end{enumerate}
	Thus the proof of the claim is complete. 

	Let $\xi_{-n}\in[x_0,y_0]$ such that $g(\xi_{-n})=\xi_0$. Then $\xi_{-n}\neq\xi_0$. We now analyze the possible relationships between $[\xi_{-n},\xi_0]$ and $[\xi_0,\xi_{n}](=g[\xi_{-n},\xi_0])$.
	\begin{enumerate}
		\item[(a)] $[\xi_{-n},\xi_0]\cap[\xi_0,\xi_{n}]=\{\xi_0\}$. 
		By Lemma \ref{lem:non_wandering_arc}, there is a minimal $k_0\geq 1$ such that $g^{k_0}([\xi_0,\xi_{n}])$ meets $[\xi_{-n},\xi_0]$. This implies that the set
		$[\xi_{-n},\xi_0]\cup\cdots\cup [\xi_{k_0n},\xi_{(k_0+1)n}]$
		possesses a loop, a contradiction. 
		\item[(b)] $[\xi_{-n},\xi_0]\subseteq [\xi_0,\xi_{n}]$ or $[\xi_{-n},\xi_0]\supseteq [\xi_0,\xi_{n}]$. By the Intermediate Value Theorem, there exists a point fixed by $g$ in $[x_0,y_0]$. This contradicts the above claim.
		\item[(c)] $[\xi_{-n},\xi_0]\cap[\xi_0,\xi_{n}]=[\xi_0,\eta_0]$ for some $\eta_0\in]\xi_{-n},\xi_0[\,\cap\, ]\xi_0,\xi_{n}[$.
	\end{enumerate} 
	
To show that Case (c) cannot happen, let $\eta_{-n}=g^{-1}(\eta_0)\in]\xi_{-n},\xi_0[$. Then $\eta_{-n}\neq \eta_0$ as $g$ has no fixed points. Note that $[\eta_0,\eta_{n}]\subseteq [\eta_0,\xi_{n}]$ and $[\eta_{-n},\eta_0]\subseteq [\xi_{-n},\xi_0]$. Thus $[\eta_{-n},\eta_0]\,\cap\,[\eta_0,\eta_{n}]=\{\eta_0\}$. An
	argument similar to Case (a) yields a contradiction.
	The proof of the proposition is complete.	
\end{proof}
\begin{lem}\label{lem:periodic}
Let $U_0,\cdots,U_{p-1}$ be an attracting (or a parabolic) cycle of Fatou domains of period $p\geq 1$, such that $f(U_k)=U_{k+1}$. If two points $x,y\in \partial U_0$ have the same itinerary with respect to a critical portrait $\wt{\mc{C}}_f$, then $x=y$. 
\end{lem}
\begin{proof}
	We first assume $x\neq y$ and then prove it by contradiction. Let $x_n=f^n(x)$ and $y_n=f^n(y)$. Then $x_n\neq y_n$ for all $n\geq 0$ as $f$ is injective in each piece of $\mc{P}_f$. We write $U_n$ as $U_{n\tu{ mod }p}(=f^n(U_0))$.
	 
	 If $U_n$ is critical, then $x_n$ and $y_n$ bound a unique subarc of $\partial U_n$, say $l_n$, such that it is compactly contained in a component of $\mb{C}\setminus\mc{R}(U_n)$. In this case, we set $k_n=0$. Otherwise, $U_n$ is non-critical; let $k_n\geq 1$ be the minimal integer such that $f^{k_n}(U_n)$ is critical and we denote by $l_n$ the pullback of $l_{n+k_n}$ by the homeomorphism
	 $$f^{k_n}:\partial U_n\to \partial U_{n+k_n}.$$
	 The arc $l_n \subseteq \partial U_n$ chosen in this way is bounded by $x_n$ and $y_n$, since $U_n,\cdots, U_{n+k_n-1}$ are non-critical.
	 
	 We claim that, for all $n\geq 0$, the map $f: l_n\to l_{n+1}$ is bijective. This holds by definition when $U_n$ is non-critical. If $U_n$ is critical, let $L_{n}$ be the component of $\partial U_{n}\setminus\mc{R}(U_n)$ containing $x_n$ and $y_n$ and let $\xi:=f^{N-n}(\partial L_n)$ with $N=n+1+k_{n+1}$. Then we get a homeomorphism 
	 $$f^{N-n}:L_n\to \partial U_{N}\setminus\{\xi\}.$$
	 We have either $l_N=f^{N-n}(l_n)$ or $l_N=\partial U_N\setminus\ol{f^{N-n}(l_n)}$. If the latter happens, then as $l_n\subseteq L_n$, the arc $l_N$ would contain $\xi$, which belongs to $\partial U_N\cap \mc{R}(U_N)$. This contradicts the choice of $l_N$. Thus $l_N=f^{N-n}(l_n)$. By the definition of $l_{n+1}$, the claim holds. 

	 Since the map $f^{kp}:\partial U_0\to\partial U_0$ eventually carries $l_{0}$ onto the whole $\partial U_0$ for large enough $k$, this contradicts the above claim, which says that $f^{kp}(l_{0})=l_{kp}\neq\partial U_{0}$. The proof of the lemma is complete.	 	  
 \end{proof}

\begin{prop}\label{prop:main}
	Let $x\neq y\in J_f$ such that $\tu{itin}(x)=\tu{itin}(y)$ with respect to a critical portrait $\wt{C}_f$. Then $x$ and $y$ lie in the boundary of a Fatou domain that is eventually iterated onto a Siegel disk.
 \end{prop}
\begin{proof}
	Let $x_n=f^n(x)$ and $y_n=f^n(y)$. By Proposition \ref{prop:partition_c}, for any $n\geq 1$, we have a homeomorphism
\begin{equation}\label{eq:home}
F^n:[x_0,y_0]\to [x_n,y_n]\subseteq\, \ol{P_{s_n}}.
\end{equation}
Note that all but countably many (the ones iterated into $E_f$) points in $[x_n,y_n]\cap J_f$ have the same itinerary. By Proposition \ref{prop:semi-buried}, $[x_n,y_n]$ has no clean subarcs. Hence $[x_n,y_n]$ passes through at least one Fatou domain, say $U_0$. According to \eqref{eq:home}, one may assume that $U_n:=F^n(U_0)$ and $[z_n,z_n']:=\ol{U_n}\cap[x_n,y_n]\,(=F^n([z_0,z_0'])\,).$ 
	
We claim that $U_n$ are eventually mapped onto Siegel disks. If this is not true, by Lemma \ref{lem:periodic}, at least one of $z_n$ and $z_n'$ is iterated into $E_f$ for all $n$. One can choose a large $n$ such that
	\begin{itemize}
		\item $U_n$ is critical and periodic of period $p$;
		\item $z_n\in\partial U_n$ is fixed by $F^p$;
		\item a ray $R$ in $\mc{R}(U_n)$ supports $U_n$ at $z_n$.
	\end{itemize}
	Since $z_n$ cannot be critical, there exist two pieces $P,P'$ of $\mc{P}_f$ such that $R\subseteq \partial P\cap\partial P'$. The closure of $P$ or $P'$ contains the whole $[x_n,y_n]$. We assume $[x_n,y_n]\subseteq\ol{P}$. Then $[x_{n+kp},y_{n+kp}]$ are contained in $\ol{P}$ for each $k\geq 1$. 
	Indeed, the non-trivial arcs 
	$$[x_{n+kp},z_{n+kp}[$$ are disjoint from $\ol{U_n}$ and approach to $z_n(=z_{n+kp})$ within a complementary component of $\mc{R}(U_n)$. The fact that $R$ supports at $(U_n, z_n)$ implies that $[x_{n+kp},z_{n+kp}[\subseteq P$. 
	
	Let $l_{n+kp}$ be the arc bounded by $z_n$ and $z_{n+kp}'$ in $\partial U_n\cap \mc{R}(U_n)$. Since for each $i\geq 0$, $[x_{n+i},y_{n+i}]$ lies in the closure of a piece of $\mc{P}_f$, we have $F^p(l_{n+kp})=l_{n+(k+1)p}$. This contradicts the fact that $F^{kp}(l_n)$ will eventually cover $\partial U_n$. The claim follows.
	
	We are left to show that $[x_n,y_n]$ is formed by two closed internal rays. For otherwise, let $U_n'\neq U_n$ be Fatou domains, having non-empty intersections with $[x_n,y_n]$. For large enough $n$, by the claim, $U_n$ and $U_n'$ are Siegel disks. The centers of $U_n$ and $U_n'$ bound an arc $\gamma$ in $[x_n,y_n]$. Let $p_1$ and $p_2$ be the periods of $U_n$ and $U_n'$, respectively. The image $\gamma'=F^{p_1p_2}(\gamma)$ is still an arc, which connects the centers of $U_n$ and $U_n'$. The absence of loops implies that $\gamma=\gamma'$. This is impossible as the actions of $F^{p_1p_2}$ on $\partial U_n$ and $\partial U_n'$ are irrational rotations. The proof of the proposition is complete.	
\end{proof}

\begin{proof}[Proof of Theorem \ref{thm:main_thm}]
	If for two angles $\theta_1$ and $\theta_2$ we have $\tu{itin}(\theta_1)=\tu{itin}(\theta_2)$ with respect to a critical portrait, then the landing points $z_i$ of $R(\theta_i)$ satisfy $\tu{itin}(z_1)=\tu{itin}(z_2)$. By Proposition \ref{prop:main}, either $z_1=z_2$ or both $z_1$ and $z_2$ belong to the boundary of the same Fatou domain, which is iterated onto a Siegel disk. The proof is complete.
\end{proof}
\begin{proof}[Proof of Corollary \ref{cor:no_wandering}]
	Suppose that $C$ is a wandering continuum.
	By replaceing $C$ by some of its iterate when necessary, one may assume that for all $k\geq 0$, $f^k(C)$ are disjoint from the finite set $E_f$. Then each $f^k(C)$ is totally contained in a piece of $\mc{P}_f$. Therefore, all the points in $C$ have the same itinerary. By Proposition \ref{prop:main}, $C$ is contained in the boundary of a Fatou domain $U$ and $U$ is eventually mapped onto a Siegel disk. It follows that $C$ cannot be wandering. The proof of the corollary is complete. 
\end{proof}

\section{Monotonicity of sets of biaccessible angles}\label{application}
As an application of the main result Theorem \ref{thm:main_thm} of this paper, we prove Theorem \ref{thm:mon} regarding the monotonicity of $\tu{Acc}(f_c)$ in the family $\mc{F}$ in this section.

Let $f_c(z)=z^2+c$ for $c\in\mb{C}$. Recall that the family $\mc{F}$ of quadratic polynomials is defined as
$$\mc{F}=\{f_c\tu{ having locally connected Julia set without a Siegel disk}\}.$$

Recall also that an angle $\theta$ in the unit circle $\mathbb{T}$ is \emph{biaccessible} (with respect to $f_c$), if there exists another $\theta'$ such that both $R(\theta)$ and $R(\theta')$ land at a common point in the Julia set (of $f_c$). 
The set of all biaccessible angles is denoted by $\tu{Acc}(f_c)$.

The notion of characteristic arcs, introduced in \cite[Lemma 2.6]{Mi00}, naturally induces a partial order in $\mc{F}$.

\begin{defi}[Characteristic arc]\label{def:characteristic}
	For each $f_c\in \mc{F}$, there are two cases:
	\begin{enumerate}
		\item[(C1)] The map $f_c$ has bounded Fatou domains. There is either a parabolic or an attracting cycle of Fatou domains, whose period is assumed to be $p\geq 1$. The return map $f_c^p$ on the closure of the critical value Fatou domain $U_0$ acts as $z\mapsto z^2$ on $\ol{\mb{D}}$. Let $w_0$ be the unique point in $\partial U_0$ fixed by $f^p$. Let $R(\eta)$ and $R(\xi)$ be the rays supporting at $(U_0,w_0)$. The \emph{characteristic sector} $S_c$ is the component of $\mb{C}\setminus\ol{R(\xi)\cup R(\eta)}$ containing $U_0$; the open subarc
		\begin{equation}\label{eq:def}
		I_c:=\{\theta\in\mb{T}: R(\theta)\subseteq S_c\}
		\end{equation} 
	  of $\mb{T}$ is called the \emph{characteristic arc} of $f_c$. The trivial case that $\xi=\eta$ is equivalent to
	  $$I_c=\mb{T}^*:=\mb{T}\setminus\{0\}\Leftrightarrow J_{f_c}\tu{ is a Jordan curve}\Leftrightarrow p=1.$$
	\item[(C2)] The map $f_c$ has no bounded Fatou domain, or equivalently $c\in J_{f_c}$. If more than one ray lands at $c$, let $S_c'$ be the component of $\mb{C}\setminus\ol{\cup R(t)}$ containing zero, with $R(t)$ running over all the rays landing at $c$;
		the \emph{characteristic sector} $S_c$ is $\mb{C}\setminus \ol{S'_c}$; the \emph{characteristic arc} $I_c$ is defined by \eqref{eq:def}. If only one ray $R(\theta)$ lands at $c$, then $I_c:=\{\theta\}$.
\end{enumerate} 
For all $f_{c}$ and $f_{c'}$ in $\mc{F}$, we say $f_{c}\prec f_{c'}$ whenever $I_c\supseteq I_{c'}$. 
\end{defi}
We remark that the preimages of an angle in $\partial I_{c}$ under $m_2:\theta\mapsto 2\theta\tu{ mod }\mb{Z}$ form a critical portrait of $f_c$.

Let $f_c\in\mc{F}$ satisfying that $I_c$ is neither $\mb{T}^*$ nor a singleton. We introduce the following notations:
\begin{enumerate}
	\item[-] Let $\gamma(a,b)$ denote the open subarc of $\mb{T}$ that starts at $a$ and ends at $b$ in the anti-clockwise direction along the circle. We have $I_c=\gamma(\eta,\xi)$ by exchanging $\eta$ and $\xi$ when necessary.
	\item[-] When $R(\theta)$ and $R(\theta')$ land at a common point $z$, we denote by $L({\theta,\theta'})$ the arc $$R(\theta)\cup\{z\}\cup R(\theta').$$
	\item[-] The preimage $H_c:=f_c^{-1}(S_c)$ is called the \emph{forbidden area} of $f_c$. We have $0\in H_c$ in Case (C1) and $0\in\partial H_c$ in Case (C2). Thus $f_c:H_c\to S_c$ is a two-to-one branched covering.
	\item[-] The preimage $m_2^{-1}(I_c)$ of $I_c$ consists of two disjoint and symmetric arcs $I_c^+=\gamma(\eta^+,\xi^+)$ and $I_c^-=\gamma(\eta^-,\xi^-)$, where $m_2^{-1}(\eta)=\{\eta^+,\eta^-\}$ and $m_2^{-1}(\xi)=\{\xi^+,\xi^-\}$. Their lengths satisfy
	\begin{equation}\label{eq:length}
	|I_c^+|=|I_c^-|=|I_c|/2.
	\end{equation}
\end{enumerate}
\begin{figure}[h]
	\begin{tikzpicture}
	\node at (0,0) {\includegraphics[width=15cm]{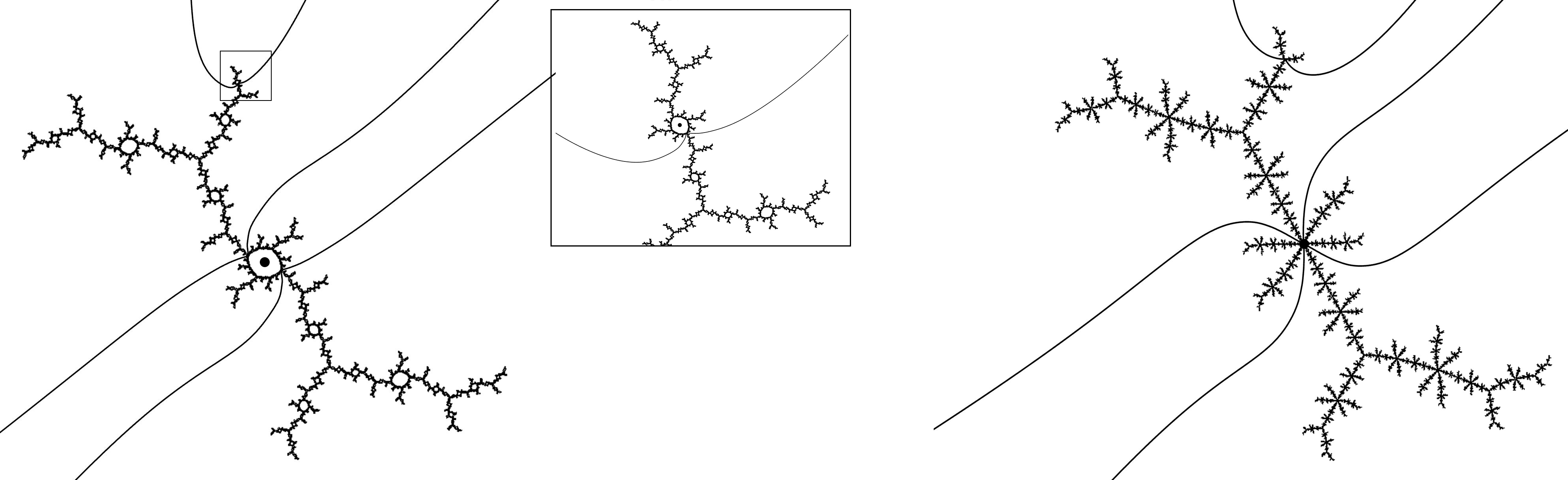}};
	\node at (3,-1.5) {$S_c^+$};
	\node at (6.75,1.25) {$S_c^-$};
	\node at (5,2.25) {$S_c$};
	\node at (5.75,2) {$R(\eta)$};
	\node at (4,2.25) {$R(\xi)$};
	\node at (3,-0.75) {$R(\eta^+)$};
	\node at (3.5,-2)  {$R(\xi^+)$};
	\node at (6.75,0.5) {$R(\eta^-)$};
	\node at (6.725,2.25){$R(\xi^-)$};
	\draw [->](-5.009,1.5)-- node [above] {} (-1.59,1.5);
	\node at (-5.25,2.25) {$S_c$};
	\node at (-6.,1.90) {$L({\eta,\xi})$};
	\node at (-6.5,-0.5) {$L({\eta^+,\xi^-})$};
	\node at (-5.5,-1.25) {$L({\xi^+,\eta^-})$};
	\node at (-3.75,0.75) {$H_c$};

	\end{tikzpicture}
	\caption{Left: Case (C1), right: Case (C2)}
	\label{fig:_illustra_properties_characteristic_arc}
\end{figure}

With all the settings above, clearly $f_{c}\prec f_{c'}$ would imply $I_{c'}\subseteq I_{c}, I_{c'}^{+}\subseteq I_{c}^+$ and $I_{c'}^-\subseteq I_{c}^-$. Before proving Theorem \ref{thm:mon}, we need the following two lemmas.
\begin{lem}
	Let $f_c\in\mc{F}$ satisfying that $I_c$ is neither $\mb{T}^*$ nor a singleton. 
	\begin{enumerate}
		\item In Case (C1), we have
			\begin{enumerate}
				\item[(1.1)] $L({\eta,\xi})$ separates $0$ and $c$;
				\item[(1.2)] $R(\eta^+)$ and $R(\xi^-)$ (resp. $R(\eta^-)$ and $R(\xi^+)$) land at the same point;
				\item[(1.3)] $H_c$ is bounded by $L({\eta^+,\xi^-})$ and $L({\eta^-,\xi^+})$. Moreover, $H_c\cap S_c=\emptyset.$
			\end{enumerate}
		\item In Case (C2), we have
			\begin{enumerate}
				\item[(2.1)] $R(\eta^+),R(\eta^-),R(\xi^+)$ and $R(\xi^-) $ land at zero;
				\item[(2.2)] let $S_c^+$ (resp. $S_c^-$) be the components of $\mb{C}\setminus L({\eta^+,\xi^+})$ (resp. $\mb{C}\setminus L({\eta^-,\xi^-})$) disjoint from $R(\xi^-)$ (resp. $R(\xi^+)$). Then $H_c=S_c^+\cup S_c^-$ and $H_c\cap S_c=\emptyset$.
			\end{enumerate}
	\end{enumerate}
\end{lem}
\begin{proof}
	(1) Let $U_0$ be the critical value Fatou domain appeared in Case (C1) of Definition \ref{def:characteristic}, and $p$ be the period of $U_0$. As $p=1$ would imply $I_c=\mb{T}^*$, we see that $p\geq 2$. Let $L_0:=L({\eta,\xi})$ and $L_k:=f^k(L_0)$, where $\eta$ and $\xi$ appeared in Case (C1) of Definition \ref{def:characteristic} satisfy $I_c=\gamma(\eta,\xi)$.
	Since the cycle $w_0,\cdots,w_{p-1}$ contains no critical point, one may let $L^+_k$ and $L^-_k$ be the two components of $\mb{C}\setminus L_k$ such that $0\in L^-_k$. 
	
	(1.1) If the statement fails, we have $\{0,c\}\subseteq L^-_0$. This would imply that all $L_{p-1}, \cdots, L_{1}$ do not separate $0$ and $c$ as well as $L_0$. Hence by Lemma \ref{lem:portrait_map}, one has
	$$(L^+_0=)L^+_{p}=\tau(L^+_{p-1})=\cdots =\tau^p (L^+_{0})$$
	and then $l(L^+_0)=2^p\,l(L^+_0),$
	which is a contradiction. Under the assumption that $c\in L_0^-$, we now show that $c\in L_k^-$ for all $k$ by induction. For $k=p-1,\cdots,1$, if $c\notin L_k^-$, then $L_k$ separates $0$ and $c$, and then $L_k\subseteq L_0^-$. It follows that
	\begin{equation}\label{eq:annular}
		L^-_{k}\subseteq L^-_0\tu{ and }L^+_0\subseteq L^+_{k}.
	\end{equation}
	On the other hand, by induction and Lemma \ref{lem:portrait_map}, the sector map $\tau$ has the orbit $$(L^-_{k}, L^+_{k})\to(L^-_{k+1}, L^+_{k+1})\to \cdots\to (L^-_p, L^+_p)(=(L^-_0, L^+_0)).$$
	Thus $l(L^+_0)=2^{p-k}l(L^+_k)$. It contradicts that $l(L^+_0)\leq l(L^+_k)$ from \eqref{eq:annular}. 
	
	$(1.2)$ Since $L_0$ is disjoint from $c$, the preimage $f_c^{-1}(L_0)$ has two components formed by the closure of four rays $R(\eta^\pm)$ and $R(\xi^\pm)$.
	By contradiction we assume that  $L(\eta^+,\xi^+)$ and $L(\eta^-,\xi^-)$ are well-defined. Let $L_{-1}^{\pm}$ be the two sectors of the portrait $\{\eta^+,\xi^+\}$ such that $L_{-1}^-$ is critical. Then $l(L^+_{-1})=|I_c^+|$ as $L_{-1}^+=\cup_{\theta\in I_c^+} R(\theta)$. Since $c\notin L_0^-$ by (1.1), we have $\tau(L_{-1}^+, L_{-1}^-)=(L_0^-,L_0^+)$. By \eqref{eq:length}, it holds that $$l(L_0^-)=2l(L_{-1}^+)=2|I_c^+|=|I_c|=l(L_{0}^+).$$
	Consequently $l(L^-_0)=\frac{1}{2}$, contradicting the fact that $l(L^-_0)>\frac{1}{2}$ as $L^-_0 $ is critical. Hence the statement follows.
	
	$(1.3)$ Let $w_{-1}$ and $\wt{w}_{-1}$ be the common landing points of the rays in $L_{-1}:=L(\eta^+,\xi^-)$ and $\wt{L}_{-1}:=L(\eta^-,\xi^+)$ respectively. Then $\{w_{-1},\wt{w}_{-1}\}=f_c^{-1}(w_0)\subseteq \partial U_{p-1}$. Since $L_0$ separates $U_0$ and $U_{p-1}$ by (1.1), the arcs $L_{-1}$ and $\wt{L}_{-1}$ are contained in the closure of $L_0^{-}$.  It is clear that $\partial H_c=L_{-1}\cup\wt{L}_{-1}$ by (1.2). To show $H_c\cap S_c=\es$, it suffices to rule out the case that $S_c\subseteq H_c$. If this case does happen, by (1.2) we have  $$\tu{either }S_c\subseteq \cup _{\theta\in I_c^+}R(\theta)\tu{ or }S_c\subseteq \cup _{\theta\in I_c^-}R(\theta).$$
	Then $|I_c|=l(S_c)\leq |I_c^+|=|I_c^-|$. This contradicts \eqref{eq:length}. The statement of (1.3) follows.
	
	$(2)$ In this case, $0\in J_{f_c}$ and the statement of (2.1) holds immediately. Since $c\in \partial S_c$, the preimage $f^{-1}_c(S_c)$ of $S_c$ has two symmetric components $S_c^+$ and $S_c^-$. Both of them are sectors based at 0 with $l(S_c^+)=l(S_c^-)=\frac{1}{2}l(S_c)$. Hence, neither of them contains $c$. According to Lemma \ref{lem:distinct_portrait}, it holds that $(S_c^+\cup S_c^-)\cap S_c=\es$. The statement of (2.2) follows.
\end{proof}

\begin{lem}\label{lem:accessible}
	Let $f_c\in\mc{F}$ satisfying that $I_c$ is neither $\mb{T}^*$ nor a singleton. Let $z_0$ be a biaccessible point whose forward orbit is disjoint from zero. Then for each sufficiently high iteration $z_n:=f^n(z_0)$ of $z_0$, we have 
	\begin{itemize}
   	\item[(1)] the forward orbit of $z_n$ is disjoint from the forbidden area $H_c$;
   	\item[(2)] any two angles whose rays land at $z_n$ have the same itinerary with respect to the partition induced by $m_2^{-1}(\nu)$ for each $\nu\in I_c$. 
	\end{itemize}
\end{lem}
\begin{proof}
	Let $R(\theta_n)$ and $R(\theta_n')$ be two distinct rays landing at $z_n$. Let $L^-_n$ and $L^+_n$ be the two components of $\mb{C}\setminus L(\theta_n,\theta_n')$ such that $0\in L_n^-$. 
	
	First, there exists a large $n$ such that $L(\theta_n,\theta_n')$ separates $0$ and $c$. For otherwise, by Lemma \ref{lem:portrait_map}, the sector map always sends $(L^-_n, L^+_n)$ to $(L^-_{n+1},L^+_{n+1})$, and then $$l(L^+_{n+k})=2^k l(L^+_n)\to\infty\tu{ as }k\to\infty,$$
	which is impossible.	
	
	Next, for such an $n$, the forward orbit $z_{n+k}$ of $z_n$ with $k\geq 0$, will never enter $H_c$. Since $H_c=f_c^{-1}(S_c)$, it suffices to show that $z_{n+k+1}$ is disjoint from $S_c$. We suppose this is not true, and assume that $z_{n+k+1}\in S_c$. Then the sector $L^+_{n+k+1}$ is compactly contained in $S_c$ and so $l(L^+_{n+k+1})<l(S_c)$. Let $n_0$ be the minimal integer in $[n+1,n+k+1]$ such that $l(L^+_{n_0})<l(S_c)$. Then by Lemma \ref{lem:portrait_map} we deduce that 
	$$\tau:(L^-_{n_0-1}, L^+_{n_0-1})\to (L^+_{n_0}, L^-_{n_0}).$$
	Thus $c\in L^+_{n_0}$ and so $S_c\subseteq L^+_{n_0}$. It implies $l(S_c)\leq l(L^+_{n_0})$, contradicting that $l(L_{n_0}^+)<l(S_c)$. Hence the statement of (1) holds. 
		
	Since $m_2^{-1}(\nu)\subseteq I_c^+\cup I_c^-$ for each $\nu\in I_c$, then the statement of (2) follows by $(1)$.
\end{proof}
\subsection{Proof of Theorem \ref{thm:mon}}
\begin{proof}[Proof of Theorem \ref{thm:mon}]
	If $I_{c}=\mb{T}^*$, then $\tu{Acc}(f_c)=\emptyset$. The theorem holds obviously. We now decompose $\tu{Acc}(f_c)$ into two disjoint subsets
	$$\tu{Acc}^+(f_c)\tu{ and }\tu{Acc}^-(f_c)$$
	such that $\theta\in \tu{Acc}^-(f_c)$ if and only if for some $n\geq 0$, the landing point of $f_c^n(R(\theta))$ is in the boundary $\partial H_c$ of $H_c$.
	
	If $I_{c'}\subsetneq I_c$, then at least one endpoint $\nu$ of $I_{c'}$ is contained in $I_c$. Thus $I_c$ is neither $\mb{T}^*$ nor a singleton. Note that $m_2^{-1}(\nu)$ is a critical portrait of $f_{c'}$. Given an arbitrary $\theta\in \tu{Acc}(f_c)$, let $z_0$ be the landing point of $R(\theta)$. If $\theta\in\tu{Acc}^+(f_c)$, by Lemma \ref{lem:accessible} (2) and Theorem \ref{thm:main_thm}, for each large $n$, all the angles whose rays landing at $z_n$ belong to $\tu{Acc}(f_{c'})$. Since
	\begin{equation}\label{eq:invariant}
	m_2^{-1}(\tu{Acc}(f_{c'}))\subseteq \tu{Acc}(f_{c'}),
	\end{equation}
	we have $\theta\in \tu{Acc}(f_{c'})$. Otherwise, $\theta\in\tu{Acc}^-(f_c)$, then for all $n$, the point $z_n$ is biaccessible. Moreover, when $n_0$ is large, the forward orbit of $z_{n_0}$ is disjoint from zero. Then one can use similar arguments as above on $z_{n_0}, \theta_{n_0}$ to show that $\theta_{n_0}\in\tu{Acc}(f_c)$. Hence $\theta\in \tu{Acc}(f_c)$ by \eqref{eq:invariant}.
	
	
	If $I_{c}=I_{c'}$, then $f_c$ and $f_{c'}$ have the same critical portraits. It holds that $\tu{Acc}^+(f_c)=\tu{Acc}^+(f_{c'})$ by Theorem \ref{thm:main_thm}.
	To show that $\tu{Acc}^-(f_c)=\tu{Acc}^-(f_{c})$, there are three cases to discuss:
	\begin{enumerate}
		\item If an endpoint of $I_c$ is periodic, then both $f_c$ and $f_{c'}$ are in Case $(C1)$. We have 
		\begin{equation}\label{eq:angles}
			\tu{Acc}^{-}(f_c)=\cup_{k\geq 1}m_2^{-k}(\partial I_c)=\cup_{k\geq 1}m_2^{-k}(\partial I_{c'})=\tu{Acc}^{-}(f_{c'}).
		\end{equation}
		\item If $\# I_c=\#I_{c'}=1$, then only one ray terminates at $c$ and $c'$ for both $f_c$ and $f_{c'}$. In this case, $\tu{Acc}^-(f_c)$ and $\tu{Acc}^-(f_{c'})$ are the iterated preimages of $\partial I_c$ and thus \eqref{eq:angles} holds again. 
		\item Otherwise, the critical values of $f_c$ and $f_{c'}$ are biaccessible. Since $f^k_c(c)$ (resp. $f_{c'}^k(c')$) are disjoint from $\partial H_c$ (resp. $\partial H_{c'}$) for all $k\geq 0$, all the angles whose rays landing at $c$ (resp. $c'$) belong to $\tu{Acc}^+(f_c)$ (resp. $\tu{Acc}^+(f_{c'})$). Thus \eqref{eq:angles} follows in this situation.
	\end{enumerate}
We complete the proof of the theorem.
	
\end{proof}

\end{document}